\documentclass[]{interact}
\usepackage{amsmath}
\usepackage[long,,nocomma]{optidef}
\usepackage[super]{nth}
\usepackage{float}
\usepackage{multirow}

\usepackage{siunitx}
\usepackage[colorlinks=false,allcolors=blue]{hyperref}

\usepackage{url}
\newcommand{\ProbOpt}[1]
{
	$\mathcal{P}_\text{#1}$
}
\usepackage[linesnumbered,ruled,vlined]{algorithm2e}

\usepackage{bm}

\usepackage{epstopdf}
\usepackage{subfigure}

\usepackage{natbib}
\bibpunct[, ]{(}{)}{;}{a}{}{,}
\renewcommand\bibfont{\fontsize{10}{12}\selectfont}

\begin{document}

 \articletype{Author’s Original Manuscript (AOM)}

\title{A Branch and Bound Based on NSGAII Algorithm for Multi-Objective Mixed Integer Non Linear Optimization Problems}

\author{
\name{A. Jaber \textsuperscript{a,b}\thanks{CONTACT A. Jaber. Email: ahmed.jaber@utt.fr}, P. Lafon\textsuperscript{a} and R. Younes\textsuperscript{b}}
\affil{\textsuperscript{a} University of Technology of Troyes, France; \textsuperscript{b} Lebanese University, Lebanon}}

\maketitle
\begin{abstract}
Multi-Objective Mixed Integer Non-Linear Programming problems (MO-MINLPs) appear in several real-world applications, especially in the mechanical engineering field. To determine a good approximated Pareto front for this type of problems, we propose a general hybrid approach based on a Multi Criteria Branch-and-Bound (MCBB) and Non-dominated Sorting Genetic Algorithm 2 (NSGAII). We present a computational experiment based on statistical assessment to compare the performance of the proposed algorithm (BnB-NSGAII) with NSGAII using well-known metrics from literature. We propose a new metric, Investment Ratio (IR), that relate the quality of solution to the consumed effort. We consider five real world mechanical engineering problems and two mathematical ones to be used as test problems in this experiment. Experimental results indicate that BnB-NSGAII could be a competitive alternative for solving MO-MINLPs.
\end{abstract}
\begin{keywords}
Multi-criteria; branch-and-bound; NSGAII; MINLP; mixed variables
\end{keywords}
\section{Introduction}
Mixed Integer Non-Linear Programming (MINLP) problems appear in several real-world applications \citep{cacchiani_branch-and-bound_2017}. 
Those are characterized by a combination of continuous variables with integer, binary, and/or categorical variables, the latter ones, being discrete, must take their values from a finite predefined list or set of categories. Categorical variables may be integers that represent the index of non-numeric values in a catalogue, such as colour, shape, or type of material \citep{abramson_mesh_2009}. A plethora of optimization methods were developed to solve these kinds of problems. They can be classified into two categories: mono-objective that aims to solve optimization problems having one objective function resulting in an optimal solution, and multi-objective methods that aim to solve multi-objective MINLPs (MO-MINLP) resulting in a set of solutions called Pareto Front.

In another manner, optimization methods can be classified into exact, approximate, or hybrid methods \citep{puchinger_combining_2005}. Exact methods are deterministic algorithms that converge towards a solution that satisfies some optimality conditions. These algorithms have been successfully applied to many engineering design problems \citep{rao_mechanical_2012}. Exact methods guarantee to find the real optimum of the problem. However, the computation cost increases prohibitively when the combinatorial space resulted from integer variables is large. Another aspect related to exact methods that solve MINLPs is that most of them rely on the concept of integer relaxation, which may be not desired in some cases. The other hand, approximate methods are stochastic in nature with probabilistic transition rules. These techniques were developed to overcome the shortcomings of exact methods. They can solve complex optimization problems resulting in a good quality solution in a relatively short time.\par
Non-Dominated Sorting Genetic Algorithm 2 (NSGAII) \citep{deb_fast_2002} is a well-known metaheuristic method that can natively handle MO-MINLPs. Any metaheuristic consists of two forces, exploration force and exploitation force. Exploration is the process of visiting entirely new regions of a search space,
whilst exploitation is the process of visiting those regions of a search space within the neighbourhood of previously visited points. Enhancing exploration force of NSGAII has been a point of interest for many researchers. A literature review was done for distinct trials for this enhancement in \citep{crepinsek_exploration_2013}. NSGAII exploration and exploitation forces depend on the cross-over and mutation mechanisms. Most researchers address those mechanisms to enhance the exploration strength.\par 
On the other hand, quite an impressive number of algorithms were reported that do not purely follow the paradigm of a single traditional metaheuristic \citep{blum_hybrid_2011}. On the contrary, they combine various algorithmic components, often originating from algorithms of other research areas on optimization. These approaches are commonly referred to as hybrid metaheuristics .\par    
The hybridization of metaheuristics with tree search techniques is probably one of the most popular lines of combining different algorithms for optimization. One of the basic ingredients of an optimization technique is a mechanism for exploring the search space, that is, the space of valid solutions to the considered optimization problem. The hybridization of metaheuristics with branch and bound concepts is rather recent \citep{blum_hybridizations_2008}. Branch and bound is a well-known generic method for exploring the domain space of mono-objective MINLPs based on the “divide to conquer” idea. It consists in an implicit enumeration principle, viewed as a search tree. The partitioning of this tree is called branching. A branch and bound algorithm produces for each subspace in the search tree an upper bound as well as a lower bound of the objective function. These bounds are used to decide if the whole subspace can be discarded, or if it has to be further partitioned.\par 
This article aims to develop and investigate the performance of a new hybrid algorithm to improve the quality of optimal solutions, i.e. Pareto front, for real mechanical engineering design optimization MO-MINLP problem by hybridization of branch and bound algorithm and metaheuristic algorithm.
\subsection{Combining mono-objective branch and bound with metaheuristics}
 A first attempt to combine branch and bound with an evolutionary algorithm has been carried in \citep{nagar_combined_1996}. At the first stage of the algorithm, branch and bound is executed up to a fixed tree-level k. Lower bounds are evaluated using a heuristic method.  The second stage consists in the execution of the evolutionary algorithm on the active nodes. In \citep{cotta_embedding_2003}, a framework based on branch and bound is presented as an operator embedded in the evolutionary algorithm. The resulting hybrid operator intelligently explores the possible children of the solutions being recombined, providing the best combination of them that can be constructed without introducing implicit mutation. In \citep{woodruff_chunking_1999}, a chunking based strategy is proposed to decide at each node of the branch and bound tree whether or not a metaheuristic is called to eventually find a better incumbent solution. In \citep{jozefowiez_bi-objective_2007}, NSGAII generates potentially Pareto optimal solutions, which are used to build sub-problems; these sub-problems are then solved using the branch-and-cut algorithm. In \citep{puchinger_solving_2004} the evolutionary algorithm is based on an order representation, specific recombination and mutation operators, and a decoding heuristic. In one variant, branch-and-bound is occasionally applied to locally optimize parts of a solution during decoding. In \citep{blum_hybridizations_2008}, branch and bound concepts are used to make the solution construction process of ant colony optimization more effective. In \citep{backer_solving_2000}, the authors only use the constraint programming system, which uses the technique of depth-first branch and bound, to check the validity of solutions and determine the values of constrained variables, not to search for solutions. Extensive literature reviews on this topic are available in \citep{puchinger_combining_2005,blum_hybrid_2011,ehrgott_hybrid_2008}.\par
\subsection{Combining multi-objective branch and bound with heuristics}
Although the branch and bound was first suggested in (1960) \citep{land_automatic_1960}, the first complete algorithm introduced as a multi-objective branch and bound, that was identified in \citep{przybylski_multi-objective_2017}, was proposed in (1983) \citep{kiziltan_algorithm_1983}. The first branch-and-bound algorithm for multi-objective mixed binary linear problems MO-MILP was proposed later in (1998) \citep{mavrotas_branch_1998}. Multi-Criteria Branch and Bound (MCBB) proposed in \citep{mavrotas_branch_1998} is a vector maximization approach capable of producing all the efficient extreme solutions of Mixed 0-1 multi-objective linear problems. For this purpose, the conventional branch and bound algorithm was properly modified, to find all the non-dominated solutions of multiple objective problems, after the exhaustive examination of all possible combinations of the binary variables.\par
In \citep{cacchiani_branch-and-bound_2017} the authors propose a heuristic algorithm based on an MCBB. It starts with a set of feasible points obtained at the root node of the enumeration tree by solving it iteratively by $\epsilon$-constraint method. Lower bounds are derived by optimally solving Non-Linear Programming problems (NLPs) by relaxing the integer constraint. Each leaf node of the enumeration tree corresponds to a convex multi-objective NLP, which is solved heuristically by varying the weights in a weighted sum approach. In \citep{florios_solving_2010}, the authors use exact as well as approximate algorithms. The exact algorithm is a properly modified version of MCBB algorithm, which is further customized by suitable heuristics. Three branching heuristics and a more general-purpose composite branching and construction heuristic are devised. In \citep{wang_heuristic_2013}, the authors develop a heuristic method based on MCBB structure. Several lower bounds, upper bounds and dominance rules are combined into the proposed heuristic method.\par
\subsection{Contributions}
In this article, we focus on multi-objective mixed-integer non- linear problems, in which the non-linear functions are non-convex ones and a subset of the variables are required to be an integer. Many real-world applications are characterized by this type of problem, especially in the mechanical engineering field. Some of these problems are required to be solved natively, i.e. the integer constraint cannot be relaxed. To determine a good approximated Pareto front for such problems, we propose a general hybrid approach based on an MCBB and NSGAII, in which to the best of our knowledge, this is the first attempt to combine MCBB with metaheuristics. The success of using metaheuristics with mono-objective branch and bound is the main motive to propose using them with the multi-criteria branch and bound. We present computational experiments based on a statistical assessment to compare the performance of the proposed algorithm with NSGAII, by solving five mechanical problems and two mathematical ones as test problems, using well-known metrics from literature in addition to a newly proposed metric that we call Investment Ratio (IR).\par
The rest of article is organized as follows. Our proposed hybrid approach is explained in section \ref{bnb_ga}. Section \ref{doe} is devoted to presenting the numerical experiment used to test the proposed algorithm. Computational results are reported and discussed in section \ref{result}. Finally, an overall conclusion is drawn in section \ref{conc}.

\section{A Branch and Bound based NSGAII}\label{bnb_ga}
In our proposed algorithm, MCBB is used to enhance the exploration force of NSGAII by investigating the mixed-integer domain space through branching it to subdomains, then NSGAII bounds each one. In this way, MCBB guides the search using the lower bounds obtained by NSGAII. We call this approach BnB-NSGAII. \par 
The general multi-Objective MINLP problem (\ProbOpt{MO-MINLP}) is defined as follows,
\begin{mini}
        {\bm{x},\bm{y}}{\bm{f}(\bm{x}, \bm{y}) = {f_1(\bm{x}, \bm{y}),\ldots, f_p(\bm{x}, \bm{y})}}{}{}
        \addConstraint{c_j(\bm{x}, \bm{y})}{\geq0,\; j=1,...,m}{}
        \addConstraint{c_j(\bm{x}, \bm{y}))}{ = 0, \; \small j=m+1,...,l}{}
        \addConstraint{l_c \leq \bm{x} \leq u_c}, \bm{x} \in  \mathbb{R}^{n_c}
        \addConstraint{l_e \leq \bm{y} \leq u_e}, \bm{y} \in \mathbb{N}^{n_e}
    \end{mini}
     In this formulation, the problem \ProbOpt{MO-MINLP} contains a total of $n = n_c + n_e$ variables, with $n_e$ integer variables and $n_c$ continuous variables. These $n$ variables are subject to take values in the domain defined by the lower  $l_c$, $l_e$ and upper $u_c$, $u_e$ bounds. It is assumed here that the functions $f_1(\bm{x},\bm{y}) . . . ,f_k(\bm{x},\bm{y}), . . . , f_p(\bm{x},\bm{y})$ and $c_j(\bm{x},\bm{y})$ are non-linear functions and not necessarily convex.\par
    \ProbOpt{MO-MINLP} is complex and expensive to solve. The general idea is thus to solve several simpler problems instead. A trivial approach is to enumerate all possible configurations of integer variables. Then for each combination, \ProbOpt{MO-MINLP} is solved as Multi-Objective continuous Non-Linear problem (\ProbOpt{MO-NLP}), which is defined as the following formulation, since all integer variables are fixed such that $\bm{y}=\bm{\bar{y}}$. 
    \begin{mini}
        {\bm{x},\bm{\bar{y}}}{\bm{f}(\bm{x}, \bm{\bar{y}}) = {f_1(\bm{x},  \bm{\bar{y}}),\ldots, f_p(\bm{x}, \bm{\bar{y}})}}{}{}
        \addConstraint{c_j(\bm{x}, \bm{\bar{y}})}{\geq0,\; j=1,...,m}{}
        \addConstraint{c_j(\bm{x}, \bm{\bar{y}}))}{ = 0, \; \small j=m+1,...,l}{}
        \addConstraint{l_c \leq \bm{x} \leq u_c}, \bm{x}\in \mathbb{R}^{n_c}
    \end{mini}
    This technique is applicable when the number of configurations is relatively small, i.e the combinatorial space of the integer variables is small. The computation time increases prohibitively when the combinatorial space is large. The growth is exponential concerning numbers of integer variables. However in \citep{mavrotas_branch_1998}, MCBB divides \ProbOpt{MO-MINLP} by constructing a combinatorial tree that aim to partition the root node problem, \ProbOpt{MO-MINLP}, into a finite number of subproblems $Pr_1,\dots,Pr_s,\dots,Pr_n$, each $Pr_s$ is considered as a node. By solving each node for each subproblem $Pr_s$ approached, one of the following scenarios must be revealed:
\begin{itemize}
        \item[$\bullet$]  $Pr_s$ is infeasible, means that no solution satisfies all constraints.
        \item[$\bullet$] There exists a candidate optimal solution, but it is useless to compute it since it can be proved that it is not better than a previously obtained one.
        \item[$\bullet$] There exists an optimal solution that must be computed.
    \end{itemize}
    In the \nth{1} and the \nth{2} cases, $Pr_s$ is discarded (pruned, fathomed) while in the \nth{3} case, a further branching of the combinatorial tree is done by dividing $Pr_s$ into farther subproblems, called children nodes. If a node cannot be divided anymore, it is called a leaf node. Leaf nodes are then solved as \ProbOpt{MO-NLP}.\par 
    
\subsection{BnB-NSGAII}
 MCBB usually uses deterministic algorithms so they can bound the node and solve the leaf. However, the computation cost increases prohibitively when the combinatorial space is large. To avoid this, we propose to use a metaheuristic method instead. The main concept of MCBB is to solve each objective function in $Pr_s$ independently, i.e. MCBB converts \ProbOpt{MO-MINLP} to mono-objective MINLP problem at each node. In contrary, BnB-NSGAII divides the domain of \ProbOpt{MO-MINLP} into smaller subdomains by constructing an exploratory tree that aims to direct the exploitation force in NSGAII to potential elements of the solution that were out of its sight. \par
 The strategic stages defining a BnB-GA algorithm are branching procedure, upper and lower bound on the optimal value of a subproblem (bounding), choice of the subproblem to investigate (fathoming) and refinement of the obtained solution.\par 
 \subsubsection{Branching}\label{branching}
Branching is a separation procedure that splits a problem domain into subsets while ensuring that the union of the subsets remains equal to the complete set. The separation principle must, therefore, allow sharing all the integer combinations in smaller sub-spaces. In \citep{cacchiani_branch-and-bound_2017},
the branching strategy is in the following manner, for each level $j$, the $j^{th}$
integer $y_j=l,\forall\; l\in\{(l_e)_j,(u_e)_j\}$, i.e. for $j=0$, at root node all integer variables are not fixed, while at the \nth{1} level, the problem is divided to $s$ subproblems, where $s$ is the number of values in $\{(l_e)_1,(u_e)_1\}$, and so on for higher levels. 
This branching strategy acts on the integer variables bounds $l_e$ and $u_e$, thus, two types of sub-spaces of combinations might result:
\begin{itemize}
        \item[$\bullet$]  Those for which $\bar{\bm{l_e}}=\bar{\bm{u_e}}$, which means that all integer variables are fixed, thus a leaf is obtained, so it is associated with a problem of the type \ProbOpt{MO-NLP}($\bm{\bar{y}}$) with  $\bm{\bar{y}}=\bm{\bar{l_e}}=\bm{\bar{u_e}}$.
        \item[$\bullet$] Those for which $\bm{\bar{l_e}}<\bm{\bar{u_e}}$, which is associated with a problem of the type \ProbOpt{MO-MINLP}($\bm{\bar{l_e}},\bm{\bar{u_e}}$) to obtain their bounds.
    \end{itemize}
 \subsubsection{Bounding and Fathoming}\label{bounding}
  Bounding a node means to evaluate its the upper and lower bounds. In a minimization problem, the lower bound of subproblem $Pr_s$ is the value that is less than or equal to every solution in $Pr_s$. While the upper bound is the best known stored feasible solution, the value of this solution is called incumbent and is stored in an incumbent list $P_N$. For a multi-objective minimization problem, upper bound for each node is considered as the set of non-dominated points $(P_N)_s$ obtained by NSGAII. While a well-known lower bound is given by the ideal point $P^I$,
 \begin{equation}
     P^I=\min f_k(\bm{x},\bm{y}); \;\;\;\;k=1,\dots,p.
 \end{equation}
 At the root node, $P_N$ is initialized by adding the initial Pareto solution obtained via NSGAII by solving \ProbOpt{MO-MINLP}. $P_N$ is then updated as NSGAII enumerates the combinatorial tree by adding the non-dominated solutions found in each children nodes.\par
 The fathoming of a subset consists in not exploring it because this one and therefore all the subsets included in it will not contain optimal Pareto front elements. For that, it suffices to compare the evaluation of this subset with the lower bounding set. If this evaluation is dominant, then this subset may be explored again or it is fathomed. A node can be fathomed if its lower bound vector is dominated by at least one of the non-dominated points of the current set $P_N$, i.e. $P^I\geq P_N$, or if the node is infeasible.\par 
 It should be noted that conventional MCBB updates the incumbent list by adding the anchor points of non-fathomed nodes. While BnB-NSGAII returns a set of Pareto solutions from each non-fathomed node. This increases the number of non-dominated solutions that the lower bound of any child node will be compared with. Hence, BnB-NSGAII provides a higher potential for convergence.\par
 The process is then iterated until all integer variables are fixed, thus a leaf is obtained. The leaf is then solved using NSGAII as \ProbOpt{MO-NLP} solver obtaining a set of solutions $(P_N)_s$ which are then added to the incumbent list $P_N$. To this end, $P_N$ may contain some elements that are dominated by $(P_N)_s$. So Pareto filtering is then applied to $P_N$ to remove dominated elements.\par 
 \subsection{Parameters Tuning}\label{bnb_nsga} 
 In NSGAII, population size, number of generations, cross-over and mutation probabilities parameters affect the performance. Tuning these parameters may result in a more desirable solution. In BnB-NSGAII, three different instances of NSGAII are called at different levels: root, node and leaf levels. This provides three levels of parameters tuning aiming for more adaptability to distinct problems. Algorithm \ref{algo} demonstrates the main loop of BnB-NSGAII, where $\mathcal{N}$ is list of stored nodes , $s$ is the current node, $s^\prime$ is child node of $s$ and $s_{max}$ is the maximum number of nodes. \par 

\begin{algorithm}[h!]
\SetAlgoLined
\DontPrintSemicolon
 Initialization\linebreak
 $P_N=\emptyset,\; \mathcal{N} \gets \emptyset$ \linebreak
 Solve \ProbOpt{MO-MINLP} by $\text{NSGAII}_\text{root}$ \linebreak
 $\mathcal{N} \gets \mathcal{N}_0$ \linebreak
 $P_N \gets (P_N)_0$ \linebreak
 $s \gets 0$ \linebreak
 $s_{max} \gets n_{nodes}-1$\;
 \While{$\mathcal{N}\neq \emptyset$ \textbf{and} $s\leq s_{max}$}{
  Choosing a node, $N_s$, in the list $\mathcal{N}$.\\
  Create a node $N_{s^\prime}$ from $N_s$ by branching\\
  
  \If{$N_{s^\prime}\neq \emptyset$}{
  
  	\eIf{$(\bar{l_e})_{s^\prime} = (\bar{u_e})_{s^\prime}$}{
  		Solve $N_{s^\prime}$ as \ProbOpt{MO-NLP}$(x_{s^\prime},\bar{y_{s^\prime}})$ by $\text{NSGAII}_\text{node}$ \;}
  	{Solve $N_{s^\prime}$ as \ProbOpt{MO-MINLP}$(x_{s^\prime},y_{s^\prime})$ by $\text{NSGAII}_\text{leaf}$ \;}
  	\If{$P^I_{s^\prime}$ dominates $P_N$}{
   $P_N \gets P_N + (P_N)_{s^\prime}$\;   
   Pareto filtering of $P_N$\\
   $\mathcal{N} \gets \mathcal{N} + \mathcal{N}_{s^\prime}$
   }
   }
   $s \gets s +1$\\
   Cleaning the list $\mathcal{N}$: clear the leaf nodes.
 }
 \caption{Main loop of BnB-NSGAII}\label{algo}
\end{algorithm}

\section{Numerical Experiment}\label{doe}
 To evaluate the performance of the proposed algorithm, we present a benchmark experiment based on a statistical assessment to compare its performance to that of NSGAII. Both are tested on  7 MO-MINLP problems from the literature.  The true Pareto solutions of these problems are known, so the Pareto fronts resulting from both algorithms can be compared to the true ones. We present evaluation criteria from literature for quantitative comparison of the performances.
\subsection{Applications}\label{app}
To demonstrate the generality of the algorithm, the presented test problems are mixed integer/binary/discrete constrained/unconstrained problems having different Pareto shapes, discrete, continuous or discontinuous. In particular, five of the selected problems are mechanical optimization problems and two of them are mathematical ones The formulations of these problems are declared in Appendix \ref{prob_form}. Although the proposed algorithm makes no assumptions on the number of objective functions, we select bi-objective problems for more comprehensive results, i.e. the Pareto front can be illustrated on a 2-dimensional objective domain.\par
\begin{figure}[h!]
    \centering
    \subfigure[Gear problem.]{%
        \resizebox*{5cm}{!}{\includegraphics{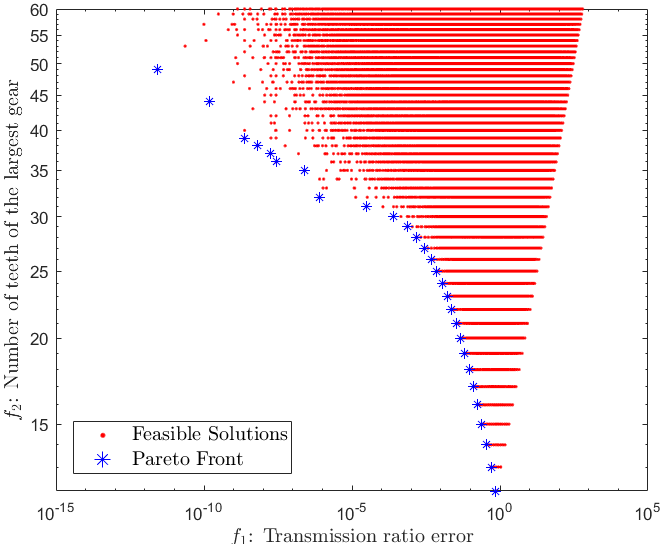}}
        \label{fig:exact_gear}
        }\hspace{5pt}
        \subfigure[Bearing problem.]{%
        \resizebox*{5cm}{!}{\includegraphics{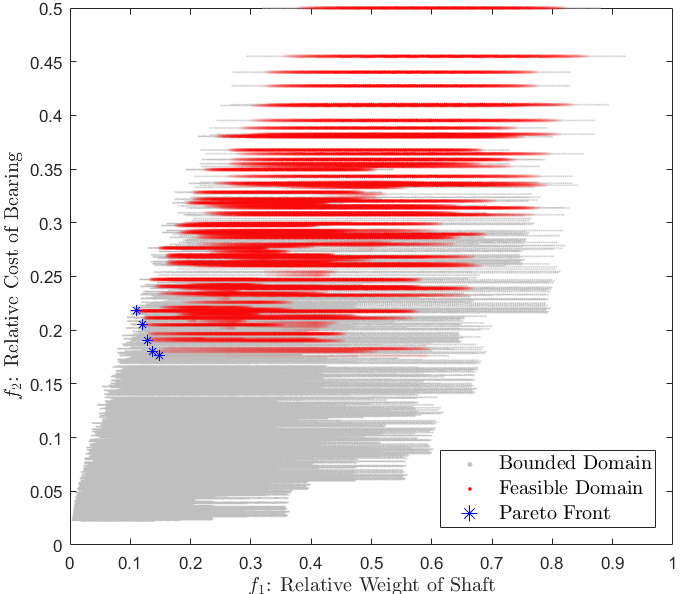}}
        \label{fig:exact_bearing}
        }
        \subfigure[Coupling problem.]{%
        \resizebox*{5cm}{!}{\includegraphics{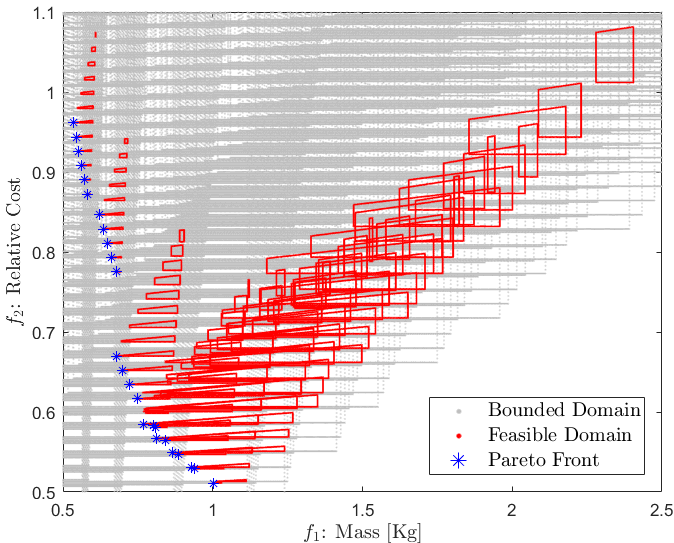}}
        \label{fig:exact_coupling}
        }\hspace{5pt}
        \subfigure[Brake problem.]{%
        \resizebox*{5cm}{!}{\includegraphics{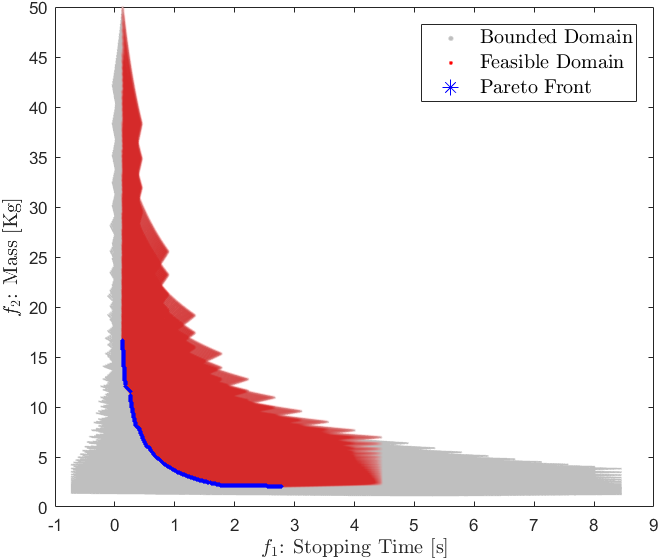}}
        \label{fig:exact_brake}
        }
    \subfigure[Truss problem.]{%
        \resizebox*{5cm}{!}{\includegraphics{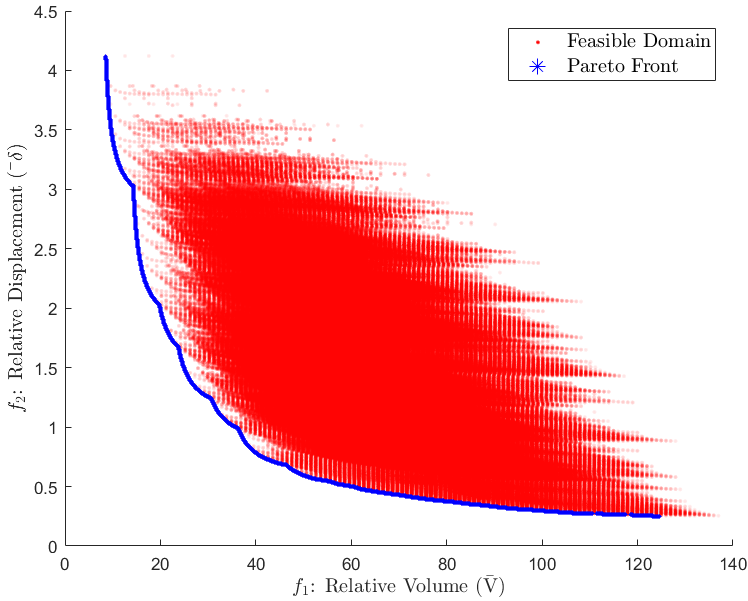}}
        \label{fig:exact_truss}
        }\hspace{5pt}
        \subfigure[Mela problem.]{%
        \resizebox*{5cm}{!}{\includegraphics{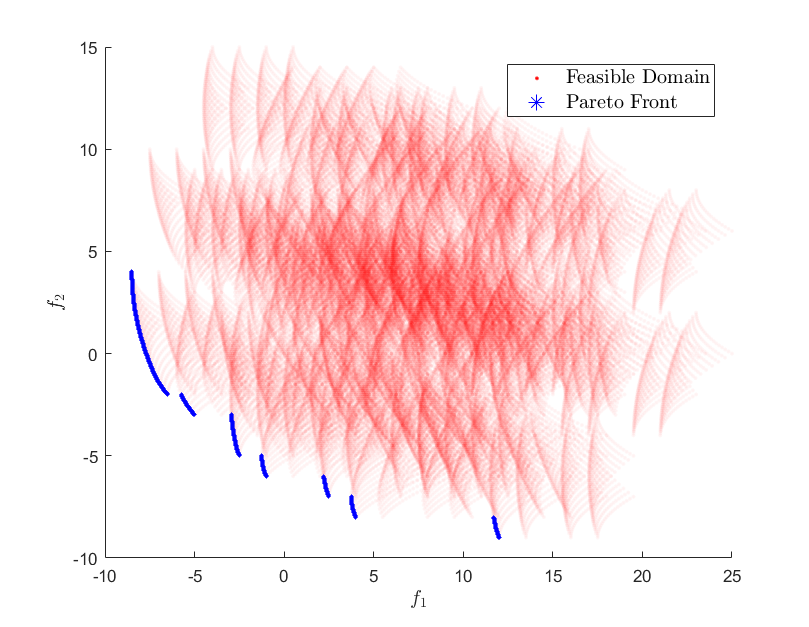}}
        \label{fig:exact_mela}
        }
    \subfigure[Tong problem.]{%
        \resizebox*{5cm}{!}{\includegraphics{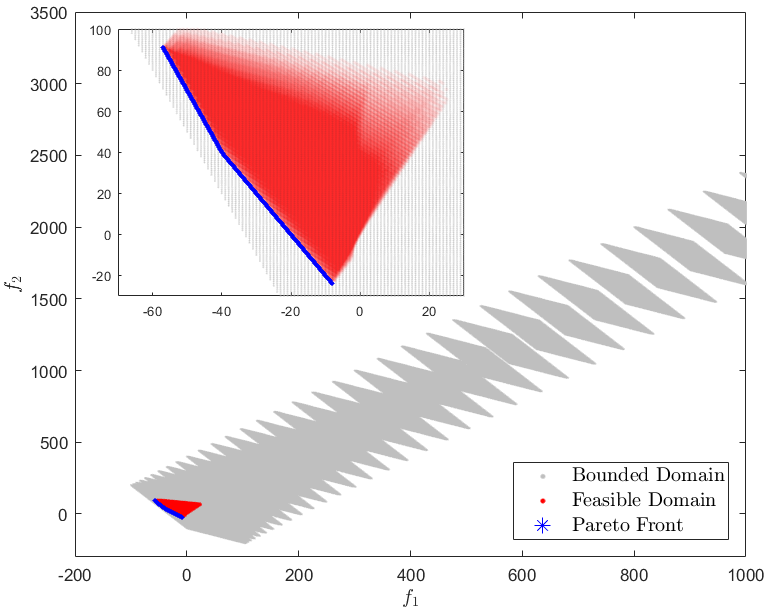}}
        \label{fig:exact_tong}
        }
    \caption{The attainable set of solutions of the test problems.}
    \label{fig:exact}
\end{figure}

 
Figure \ref{fig:exact} illustrates the exact representation of the seven problems presented in this section. For each problem, the overall domain projection on the objective space is illustrated in grey, while the feasible one in red. The true Pareto for each problem is plotted in blue. The method used to obtain the true Pareto solution is explained in Appendix \ref{true_par}. Table \ref{test_prob} summarizes various properties of the presented problems.\par 
The gear train problem proposed in \citep{deb_mechanical_2000}, represents a compound gear train which is to be designed to achieve a specific gear ratio between the driver and driven shafts. The objective of the gear train design is to find the number of teeth in each of the four gears to minimize (i) the error between the obtained gear ratio and a required gear ratio of and (ii) the maximum size of any of the four gears. The solution of this unconstrained pure integer problem is illustrated in Figure \ref{fig:exact_gear} showing that its Pareto front is discrete.\par
The ball bearing pivot link was initially proposed in \citep{giraud_optimal_1999} as a mono-objective problem, then reconstructed in \citep{el_samrout_hybridization_2019} to bi-objective one.  In this problem, the aim is to minimize the relative mass of the system composed of a shaft and two ball bearings, and relative cost of bearings. Figure \ref{fig:exact_bearing} shows the complexity of the domain. It also shows that the Pareto front of this problem is discrete.\par
The coupling with bolted rim proposed in \citep{giraud-moreau_comparison_2002} for torque transmission was reconstructed in \citep{el_samrout_hybridization_2019} to a bi-objective problem for which goal is to dimension the mechanism in a way that minimizes its cost and weight. Figure \ref{fig:exact_coupling} shows a part of the bounded domain to focus on the feasible domain that is represented as patches of solutions. The Pareto solution of this problem is a discrete one. 
In the disk brake design problem, reported in \citep{osyczka_new_1995,tong_multi-objective_2016}, there are two design objectives: minimize  the mass of the brake and minimize the stopping time. The non-convex discontinuous Pareto front is observed in figure \ref{fig:exact_brake}.\par 
Next, bi-criterion optimization of the nine bar truss was considered in \citep{mela_algorithm_2007}. The goal is to minimize the  material volume $V$ and the vertical displacement $\delta$ of the loaded node simultaneously. It is generally known that material volume and flexibility are strongly conflicting quantities in structural design. Figure \ref{fig:exact_truss} shows the minimal points of all the components of the problem. The convex minimal curves of the components intersect each other, resulting in a highly non-convex minimal curve for problem.\par 
Mela problem is a mathematical problem proposed in \citep{mela_algorithm_2007}. It is unconstrained and it is nearly linear, the quadratic term being the only nonlinear part of the problem. Figure \ref{fig:exact_mela} shows that even a mildly nonlinear  problem can have a very complicated attainable set, there are many locally minimal points. The  set of minimal points is not connected, but consists of subsets of minimal points of seven components.\par 
The analytical problem, which is adapted from \citep{dimkou_parametric_1998,tong_multi-objective_2016}, is referred to as Tong problem. In this problem, the continuous variables are unbounded. However their bound can be deduced analytically from the constraints. To examine the exploration force of the optimization method, the bounds are set to generate a large-scale domain as shown in figure \ref{fig:exact_tong}.\par 
\begin{table}[h!]
    
    \tbl{Test problems properties}
    {\begin{tabular}{lcclcl}         
            \toprule
            \multicolumn{1}{l}{\textbf{Problem}} &
            \multicolumn{1}{l}{\textbf{\# Variables}} & \multicolumn{1}{l}{\textbf{\# Integer Vars}} & \multicolumn{1}{l}{\textbf{Var Type}} & \multicolumn{1}{l}{\textbf{\# Constraints}} & \multicolumn{1}{l}{\textbf{Pareto Type}} \\ \midrule
            \textit{Gear}   & 4  & 4  & Integer  & 0  & Discrete  \\
            \textit{Bearing} &  4   & 2  & Discrete  & 10  & Discrete \\
            \textit{Coupling} &  4 & 2  & Discrete  & 4  & Discrete  \\
            \textit{Brake}    & 4  & 1 & Integer  & 5    & Discontinuous  \\
            \textit{Truss}    & 9 & 6   & Discrete  & 0  & Discontinuous\\
            \textit{Mela}    & 10  & 8  & Binary & 0  & Discontinuous \\
            \textit{Tong}    & 6  & 3  & Binary  & 9  & Continuous  \\ 
            \bottomrule
    \end{tabular}}
    \label{test_prob}
\end{table}

\subsection{Evaluation Criteria}\label{eval_crit}
The evaluation of the optimization method is a necessary task that should be conducted properly. In this section, we present three ways to evaluate a multi-objective solver; namely accuracy metrics, computational effort and robustness.\par
Accuracy metrics are a sort of performance indicators that try to assess the quality of solutions in terms of precision. Several classifications of metrics can be found in the literature including convergence-based metrics, diversity-based indicators \citep{audet_performance_2018} and others. In this article, true Pareto front based metric will be used, since the true Pareto solutions of the proposed benchmark problems are known.\par
\subsubsection{Cardinality Indicators}
These metrics focus on the number of non-dominated points generated by a given algorithm. Cardinality indicators were reviewed in \citep{audet_performance_2018}, among them is the overall non-dominated vector generation (ONVG). ONVG is the number of the Pareto front approximation generated by the algorithm:
\begin{equation}
\text{ONVG}(S)=\mid S\mid
\end{equation}
Where $\mid S\mid$ is the number of points of the approximation Pareto set. Another indicator called the purity metric is found in \citep{custodio_direct_2011}. It is used when comparing two or more approximate methods in the absence of the knowledge of the true Pareto. This indicator is slightly modified to involve the known true Pareto front. Let $\mid P\mid$ be the number of elements in the true Pareto set and $\mid F\mid=\mid S\mid \cup \mid P\mid$, and then removing from $\mid F\mid$ any dominated points. The modified purity metric is then given by:
\begin{equation}
\text{Purity}(S)=\frac{\mid{}S\mid\cap\mid{}F\mid}{\mid{}S\mid}
\end{equation}
Therefore, this value represents the ratio of non-dominated points over the solution set. This metric is thus represented by a number between zero and one. Higher values indicate a better Pareto front in terms of the percentage of non-dominated points.\par 
\subsubsection{Convergence Metrics}
Generational Distance (GD) \citep{van_veldhuizen_measuring_2000} and Inverted Generational Distance (IGD) \citep{li_multiobjective_2009} are well-known metrics with regard to convergence metrics. They compute the average distance from the approximated Pareto to the true one. For bi-objective problem, GD is defined as follows:
\begin{equation}
\text{GD} =\dfrac{\sqrt{\sum_{i=1}^{\mid{}S\mid}{d_i^2}}}{\mid{}S\mid}
\end{equation}
Where $d_i$ is the euclidean distance in objective function space between each approximated solution to the nearest true solution. A value of $GD = 0$ indicates that all the approximated solutions are on the true Pareto front. IGD substitute $\mid{}S\mid$ by $\mid P\mid$.\par
In the case of discrete Pareto, IGD indicator is used, while for continuous/discontinuous Pareto, GD indicator is used. Effect of the method used to obtain the true Pareto on this metric is discussed in Appendix \ref{true_par}.\par  
\subsubsection{Distribution and spread indicators}
The spread ($\Delta$) is another accuracy metric criteria introduced in \citep{deb_multi-objective_2001}. It includes information about both uniformity and spread. The formulation of this metric for a bi-objective problem is:
\begin{equation}
\Delta=\dfrac{d_f + d_l + \sum_{n=1}^{n-1}{|d_i-\bar{d}|}     }{d_f+d_l+(n-1)\bar{d}}
\end{equation}
where $d_f$ and $d_l$ are the euclidean distances between the extreme solutions in the non dominated set obtained using the algorithm and the real Pareto front, $n$ is the number of approximated solutions $n=\mid{}S\mid$, $d_i$ is the euclidean distance between two consecutive solutions and $\bar{d}$ is the average of all distances $d_i$, $i \in [1, n-1]$. The lower the spread indicator is, the better the solution is. Nevertheless, this indicator will never be zero, unless the spread of the true Pareto $\Delta_\text{True}$ is equal to zero and the solver succeeded to find it (see Appendix \ref{true_par}). Thus, for more comprehensible evaluation, relative spread formulated in equation \ref{eq_dD} is used,
 \begin{equation}\label{eq_dD}
    d\Delta= \mid\Delta_\text{True}-\Delta_\text{Approximate}\mid
 \end{equation}
\subsubsection{Computational Effort}
The quantification of computational effort used by the algorithm to produce the optimal Pareto front is a good way to assess its effectiveness. In \citep{talbi_metaheuristics:_2009}, the author considers the number of objective function evaluations is a good representation of the computational effort .\par 
\subsubsection{Robustness}
The robustness of the algorithm should also be measured. This is done by running it several times to compensate for the impact of random parameters, which may cause the fluctuation of solution results. The lower the variability of the obtained solutions the better the robustness \citep{montgomery_design_2017}.\par 

\subsubsection{Investment Ratio}
The necessity of more accurate solutions with the toleration of more computational effort, a new metric that measures the potential improvement of the solution concerning computational cost is defined. In our case, the quality is the accuracy of the solution, while the cost is the computational effort. In the purpose of relating the quality and the cost in a single indicator that can represent how much the quality increased with respect to the cost, we propose the investment ratio indicator. Let $Q_1$ and $Q_2$ be the quality indicators resulting from 2 solvers being compared, and $C_1$, $C_2$ be their cost. Let $q=Q_2/Q_1$ be quality ratio and $c=C_2/C_1$ be the computational cost ratio. The proposed Investment Ratio indicator (IR) is designed to regard various cases and is defined in equation \ref{eq_IR}.
\begin{equation}\label{eq_IR}
\text{IR}(q,c)=\left\{\begin{matrix}
\frac{q}{c} & \text{if} ~ ~ q\geq1\\ 
-\frac{c}{q} & \text{if} ~ ~ q<1
\end{matrix}\right.
\end{equation}
The positive value of $\text{IR}$ indicates that the quality is enhanced regardless of the cost ratio and vice versa. $\text{IR}\geq 1$ indicates a good investment, as $Q_2>Q_1$ by improvement ratio higher than the cost one. For $\text{IR}\in]0,1[$, this means that although the cost ratio is higher than the quality ratio, still the quality is enhanced. Here, the investment cannot be assessed obviously since it depends on the weights of the desire of the investor for a certain problem. $\text{IR}=-1$ indicates that the quality is reduced by the same cost reduction ratio, which is also an acceptable investment. Next, $\text{IR}\in]-1,0[$  means that solver 2 obtained lower quality with lower computational effort, yet the cost is reduced by a higher ratio than the quality degradation ratio, this value also might be acceptable depends on the investor desire. Finally, $\text{IR}<-1$ indicates a bad investment in the absence of profitable improvement in either cost or quality.\par 
In this experiment, GD and $\Delta$ are used as quality indicators, but since equation \ref{eq_IR} assumes that the higher values of $q$ mean higher quality, $q$ is defined as the geometric mean between the inverse of the spread and GD as stated in equation \ref{eq_q}. In particular, the mean value of the metrics is used,
 \begin{equation}\label{eq_q}
 q=\sqrt{\left(\frac{(GD_\text{NSGA})_\text{mean}}{(GD_\text{BnB})_\text{mean}}.\frac{(\Delta_\text{NSGA})_\text{mean}}{(\Delta_\text{BnB})_\text{mean}}\right)}
 \end{equation}
\subsection{Benchmark Experiment}
Both NSGAII and BnB-NSGAII performances are affected by selecting different values of tuning parameters. In the purpose of impartial comparison between them, different tuning parameters are applied to each method on each of the 7 problems independently. To examine the robustness of the results, each problem is solved 20 times by each method for each combination. \par
In this experiment, the effect of tuning the parameters on NSGAII and BnB-NSGAII are tested considerably, mainly for population size and the number of generations for their direct effect on both quality of the solution and the computational effort consumed. It should be noted that the number of generations is affected by two parameters. First, the maximum allowable generations which terminate NSGAII process even if the optimality condition is not attained. Second, the process is terminated if the optimality condition didn't change over $m$ generations, $m$ being the maximum stall generations value to be tuned. This test is done by varying these parameters iteratively seeking for the best performance of NSGAII and BnB-NSGAII. The test considers the three tuning levels of BnB-NSGAII mentioned previously. 
 Table \ref{table_param} shows the values of distinct parameters used for both algorithms.\par 
Best to best comparison is commonly used to compare the performances of distinct solvers. However, the best solution may be an outlier. For that, the presented benchmark is based on the distribution of values for each metric over the 20 iterations.\par
\begin{table}[H]
	\tbl{Parameters used for NSGAII and BnB-NSGAII algorithms}{
	\begin{tabular}{ll}
		\hline
		\textbf{Parameters}               & \multicolumn{1}{l}{\textbf{Value}}  \\ \hline
		Cross over probability             & 0.9                                 \\
		Mutation Probability              & 0.95                                \\
		Population size                   & variable: 20 $\rightarrow$ 1000               \\
		Allowable generations     & variable: 20 $\rightarrow$ 1000               \\
		Stall generations         & variable: 20 $\rightarrow$  Allowable   generations \\
		Constraint handling     & Legacy method\textsuperscript{a}                     \\ 
		Crossover operator				& Simulated Binary crossover (SBX)\textsuperscript{b}   \\
		ETAC & 100
		\\
		 Mutation operator				&  Partially-mapped crossover (PMX)\textsuperscript{b} \\
		
		ETAM & 10
		\\
		 \hline
	\end{tabular}}
	\tabnote{\textsuperscript{a}\citep{deb_fast_2002} \\
	\textsuperscript{b}\citep{maruyama_parametric_2017}.}
\label{table_param}
\end{table}

\section{Results and discussion}\label{result}
\begin{figure}[h!]
	\centering
	\includegraphics[width=\linewidth,height=7cm]{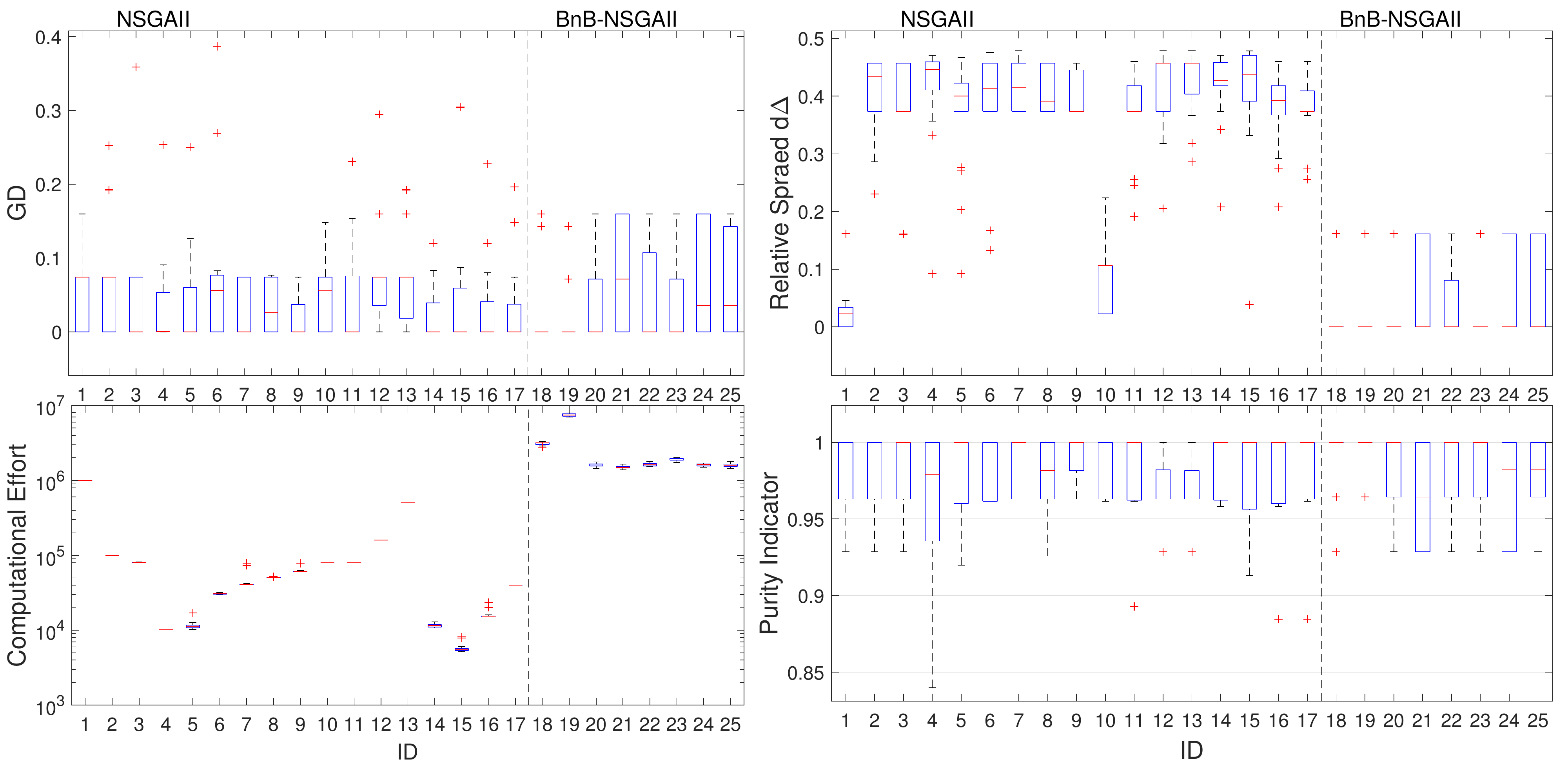}
	\caption{Distribution of metrics for gear problem.}
	\label{fig:results_gear}
\end{figure}

\begin{figure}[h!]
	\centering
	\includegraphics[width=\linewidth,height=7cm]{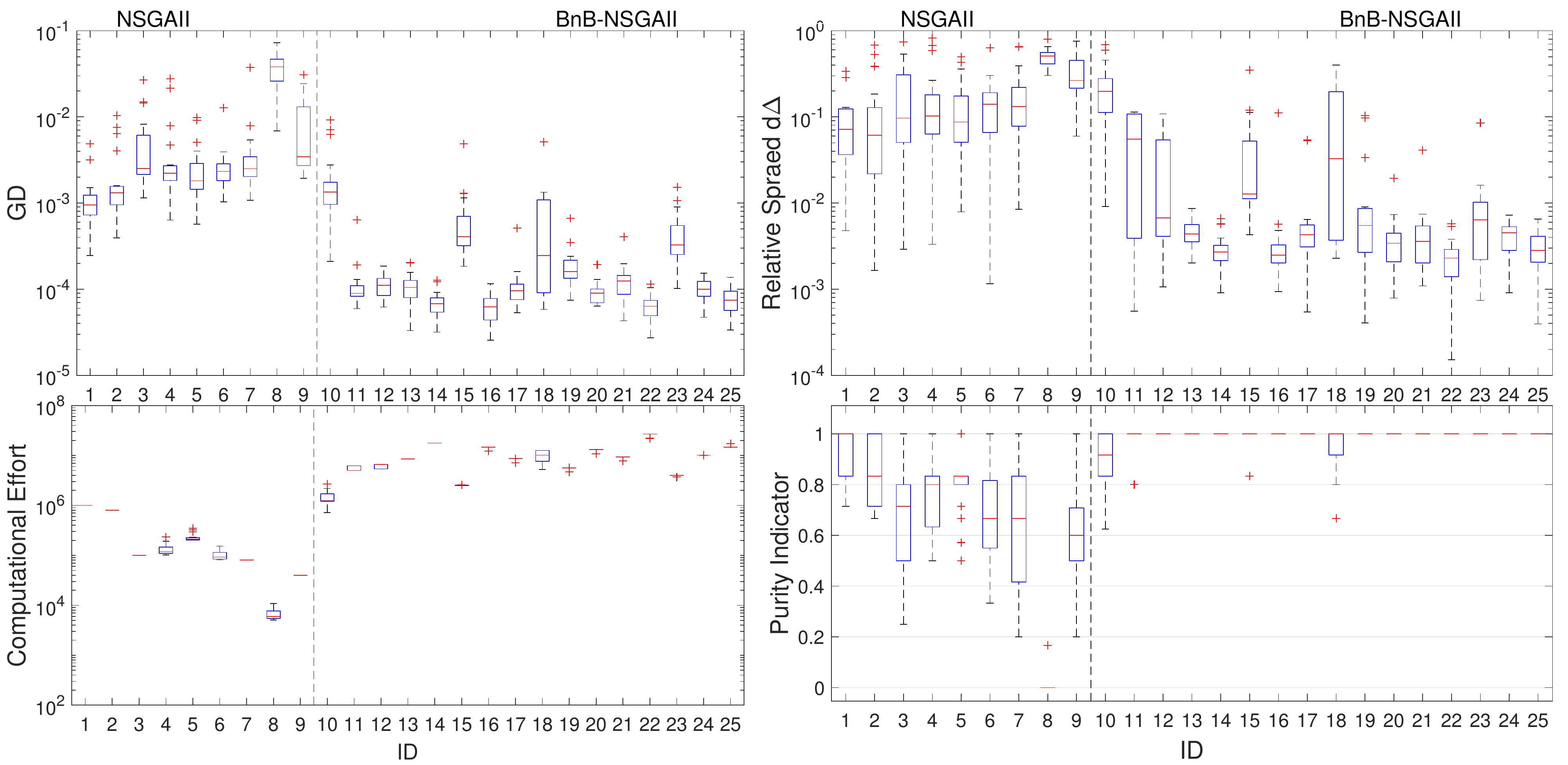}
	\caption{Distribution of metrics for bearing problem.}
	\label{fig:results_bearing}
\end{figure}

\begin{figure}[h!]
	\centering
	\includegraphics[width=\linewidth,height=7cm]{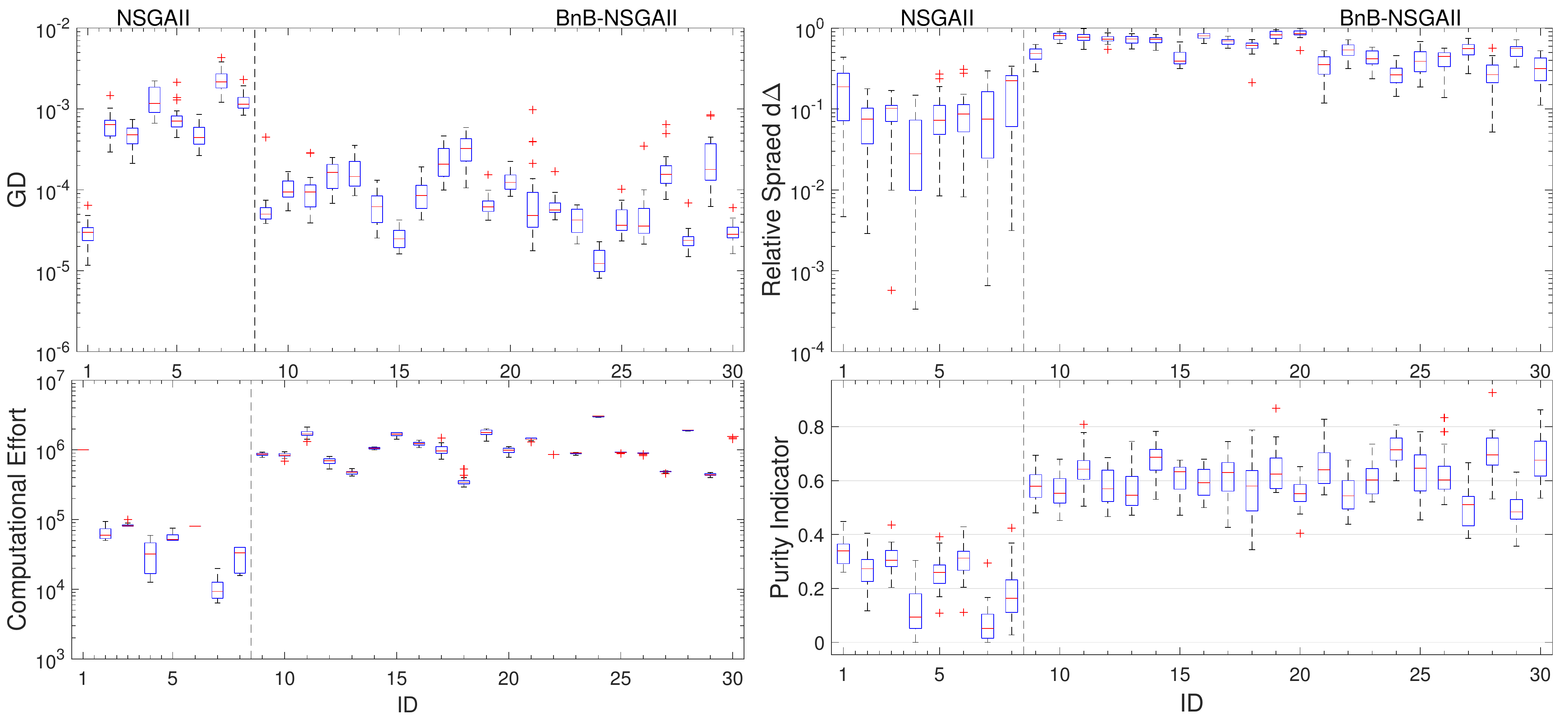}
	\caption{Distribution of metrics for coupling problem.}
	\label{fig:results_coupling}
\end{figure}

\begin{figure}[h!]
	\centering
	\includegraphics[width=\linewidth,height=7cm]{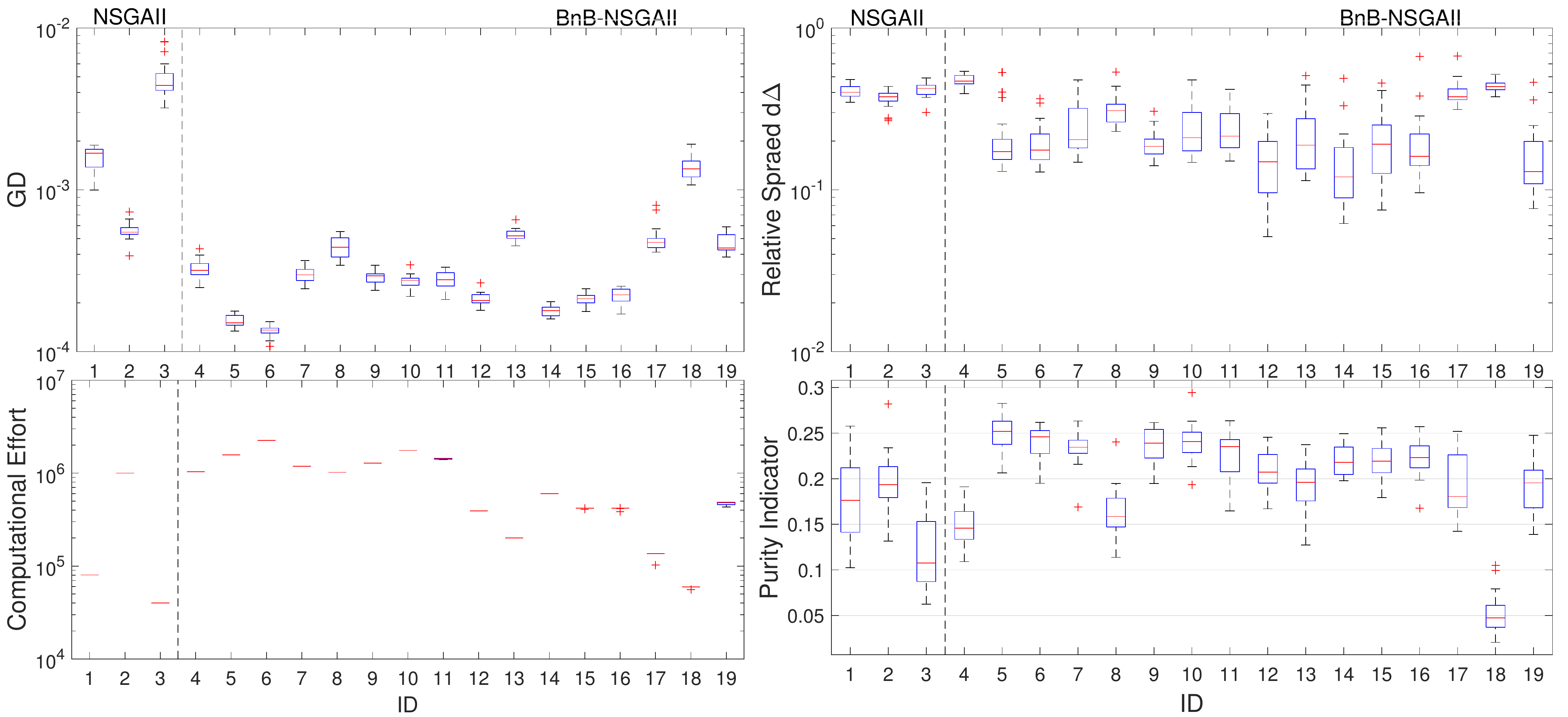}
	\caption{Distribution of metrics for disk brake problem.}
	\label{fig:results_brake}
\end{figure}

\begin{figure}[h!]
	\centering
	\includegraphics[width=\linewidth,height=7cm]{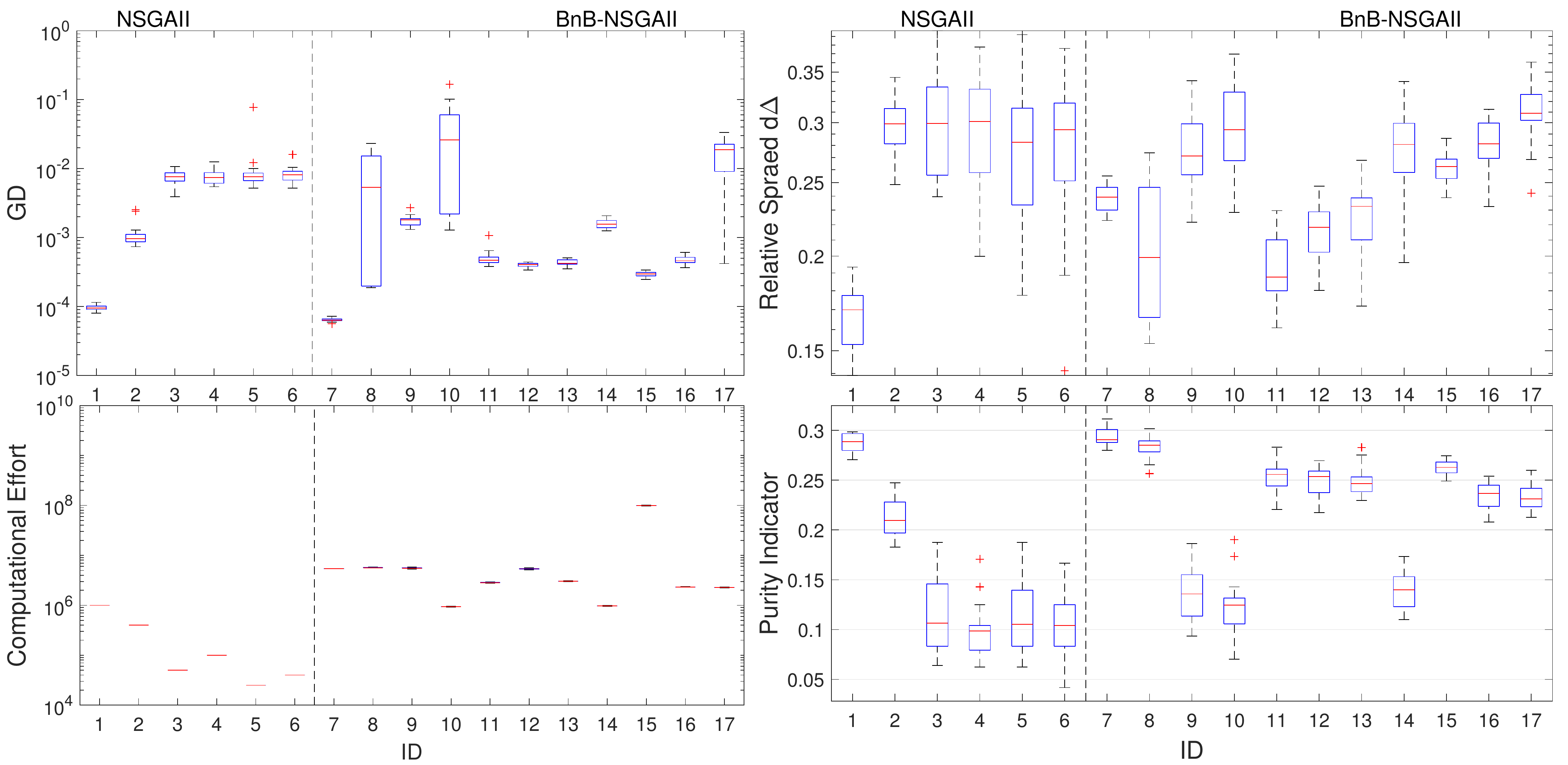}
	\caption{Distribution of metrics for truss problem.}
	\label{fig:results_truss}
\end{figure}

\begin{figure}[h!]
	\centering
	\includegraphics[width=\linewidth,height=7cm]{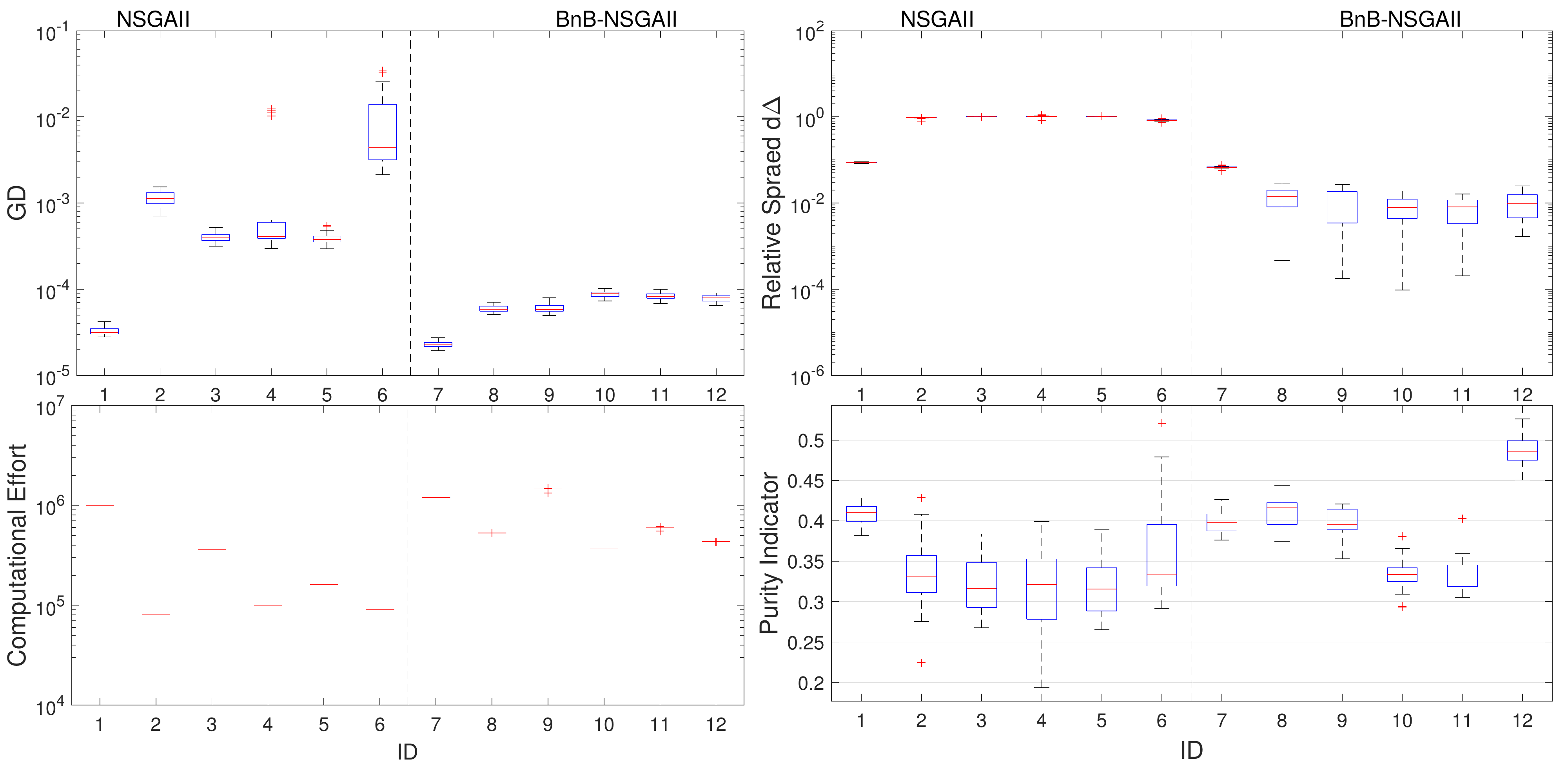}
	\caption{Distribution of metrics for Mela problem.}
	\label{fig:results_mela}
\end{figure}

\begin{figure}[h!]
	\centering
	\includegraphics[width=\linewidth,height=7cm]{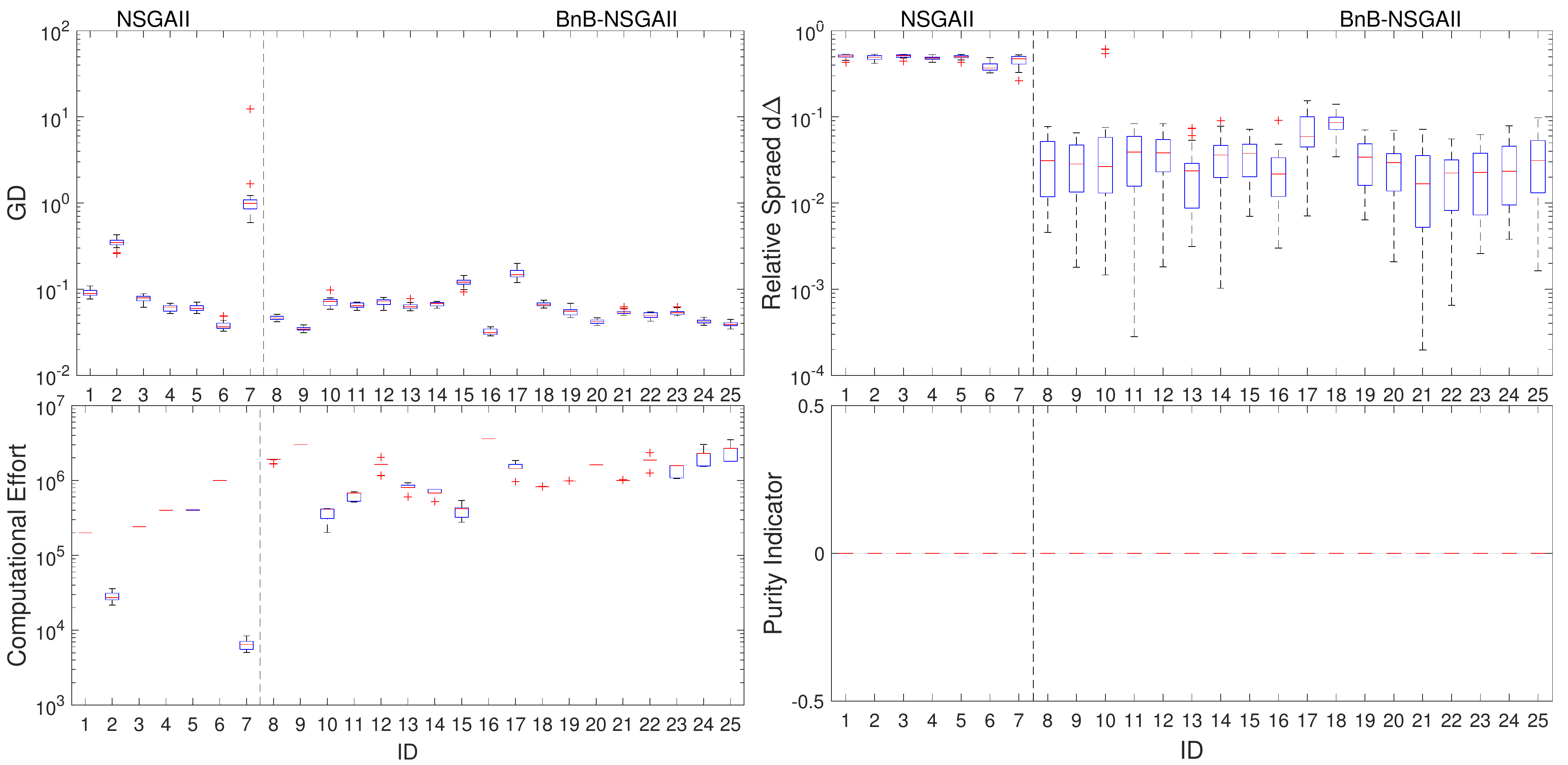}
	\caption{Distribution of metrics for Tong problem.}
	\label{fig:results_tong}
\end{figure}

Figures \ref{fig:results_gear} to \ref{fig:results_tong} show the distribution of GD, Spread, Purity, and the number of evaluations over the 20 iterations for each combination of parameters, the latter being represented as parameter ID.  It should be noted that an outlier is detected when a value is more than 1.5 times the interquartile range away from the top or bottom of the box.\par 
As expected, Figures \ref{fig:results_gear} to \ref{fig:results_tong} show that parameters tuning has significant effect on the performance of NSGAII and BnB-NSGAII. Hence, in the purpose of comprehensive discussion, best performance of NSGAII in terms of solution convergence and uniformity is selected for each problem. While for BnB-NSGAII, parameters resulting in 1$)$ most accurate results and 2$)$  best compromise between solution quality and computational effort (i.e best investment) are selected to be discussed. Tables \ref{nsga_param} and \ref{bnb_nsga_param} show the selected parameters for NSGAII and BnB-NSGAII respectively.\par 
\begin{table}[h!]
\tbl{Selected parameters for NSGAII}{
\begin{tabular}{llccc}
\hline
\textbf{Problem} & \textbf{ID} & \multicolumn{1}{l}{\textbf{Population size}} & \multicolumn{1}{l}{\textbf{Allowable generations}} & \multicolumn{1}{l}{\textbf{Stall generations}} \\ \hline
Gear             & 1           & 1000                                         & 1000                                               & 1000                                            \\
Bearing          & 1           & 1000                                         & 1000                                               & 1000                                            \\
Coupling         & 1           & 1000                                         & 1000                                               & 1000                                            \\
Brake            & 2           & 1000                                         & 1000                                               & 1000                                            \\
Truss            & 1           & 1000                                         & 1000                                               & 1000                                            \\
Mela             & 1           & 1000                                         & 1000                                               & 1000                                            \\
Tong             & 6           & 1000                                         & 1000                                               & 1000                                            \\ \hline
\end{tabular}}
\label{nsga_param}
\end{table}

\begin{table}[h!]
\tbl{Selected parameters for BnB-NSGAII}{
\begin{tabular}{lcccccccccc}
\hline
\multicolumn{1}{c}{\multirow{2}{*}{\textbf{Problem}}} & \multirow{2}{*}{\textbf{ID}} & \multicolumn{3}{c}{\textbf{Root}}                               & \multicolumn{3}{c}{\textbf{Node}}                               & \multicolumn{3}{c}{\textbf{Leaf}}                               \\ \cline{3-11} 
\multicolumn{1}{c}{}                                  &                              & \textit{Pop} & \textit{Gen} & \textit{Stall} & \textit{Pop} & \textit{Gen} & \textit{Stall} & \textit{Pop} & \textit{Gen} & \textit{Stall} \\ \hline
\multirow{2}{*}{Gear}                                 & 20                            & 100&800&500&50&100&20&50&20&20\\               & 18
&100&800&500&100&50&20&50&20&20\\
 \cline{2-11} 
\multirow{2}{*}{Bearing}                              & 19 
&50&800&800&50&800&300&50&800&300\\                   & 22
&50&800&800&50&800&800&100&800&800\\
 \cline{2-11} 
\multirow{2}{*}{Coupling}                             & 25                           &10&100&100&10&150&100&10&800&200\\               & 24
&10&100&100&10&100&100&10&800&800\\
 \cline{2-11} 
\multirow{2}{*}{Brake}                                & 12                          &50&200&200&50&200&200&100&200&200\\              & 6
&1000&1000&1000&50&20&20&100&800&500\\
 \cline{2-11} 
\multirow{1}{*}{Truss}                                & 7                        &1000&1000&1000&100&20&20&100&20&20\\
 \cline{2-11} 
\multirow{2}{*}{Mela}                                 & 10                      &200&800&100&50&20&20&50&200&200\\                    & 7
&1000&1000&1000&50&20&20&50&200&200\\
 \cline{2-11} 
\multirow{2}{*}{Tong}                                 & 10                            &200&1000&1000&50&20&20&100&800&500\\            & 16
&200&1000&1000&200&1000&1000&200&1000&1000\\
 \hline
\end{tabular}}
\label{bnb_nsga_param}
\end{table}
Compared to NSGAII, BnB-NSGAII provides better performance in terms of GD and spread indicators for all the problems except for truss problem which has higher spread. Regarding purity indicator, BnB-NSGAII overcomes NSGAII for all the problems except for Mela problem which has slightly lower ratio. For gear problem, Figure \ref{fig:results_gear} shows that BnB-NSGAII, $\text{ID}=18$, is able to find the exact solution of this problem robustly. It is also observed that BnB-NSGAII succeeds to find a better solution with lower computational effort for coupling, brake, Mela and Tong problems in Figures \ref{fig:results_coupling}, \ref{fig:results_brake}, \ref{fig:results_mela}, \ref{fig:results_tong} for $\text{ID}=25,12,10,10$ respectively.\par
To specifically explore the computational efficiency of the BnB-NSGAII algorithm, investment ratio is calculated for each of the selected combinations . Table \ref{tab_results} shows the mean values of GD, relative spread d$\Delta$, number of solutions in the approximated Pareto front $\mid S\mid$, purity ratio and computational effort after removing the outlier points, in addition to the investment ratio indicator. It is observed that BnB-NSGAII is able to provide a "good investment" ($\text{IR}>1$) for all problems except for truss problem. Therefore, the effectiveness of BnB-NSGAII is justified for 6 out of 7 problems. $\text{IR}=\infty$ in gear problem means that BnB-NSGAII is able to find a solution with best quality, i.e. the true solution. \par
\begin{figure}[H]
	\centering
	\subfigure[Gear problem.]{\resizebox{5cm}{!}{\includegraphics{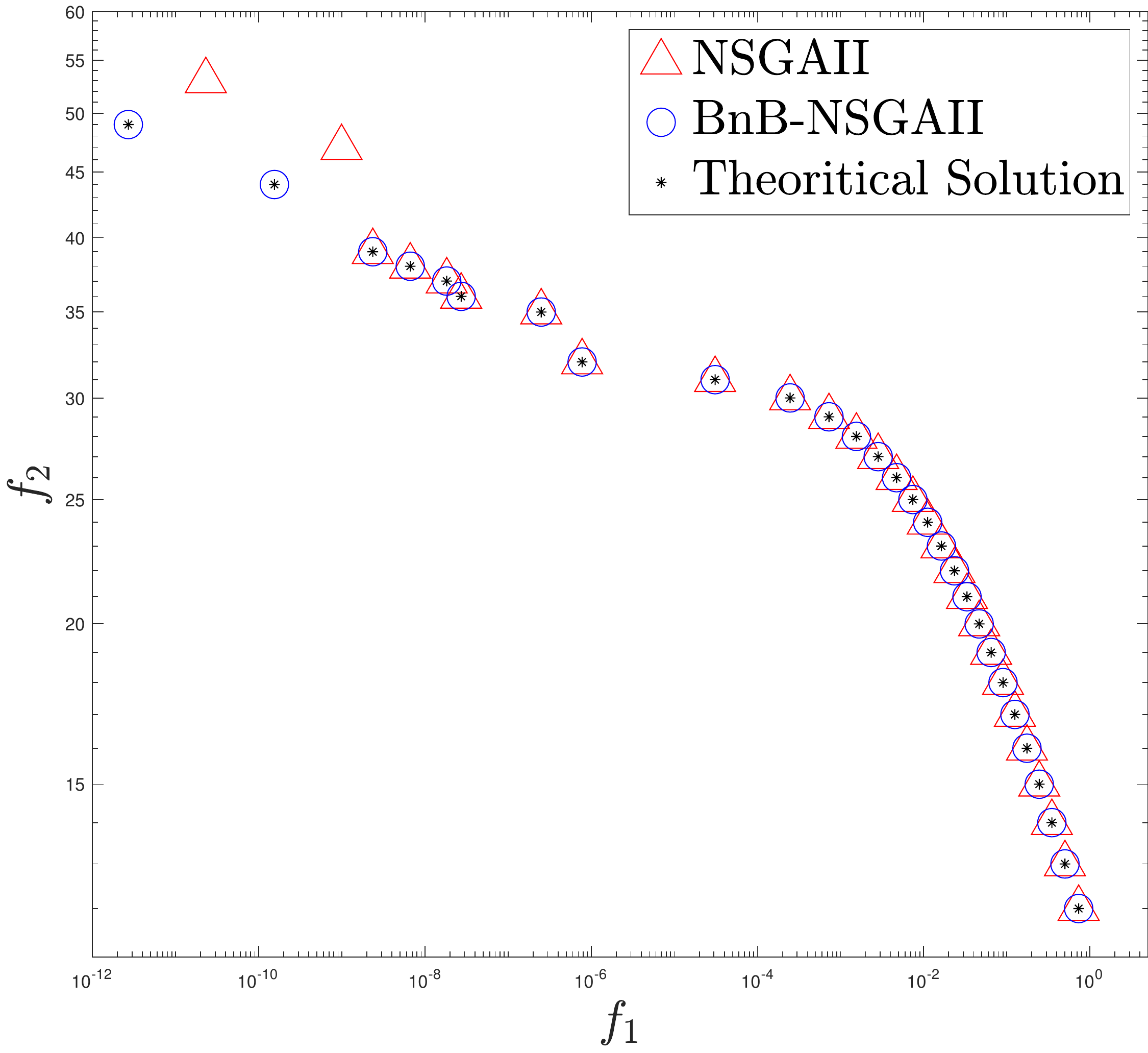}}
	\label{fig:pareto_gear}
	}\hspace{5pt}
	\subfigure[Bearing problem.]{\resizebox{5cm}{!}{\includegraphics{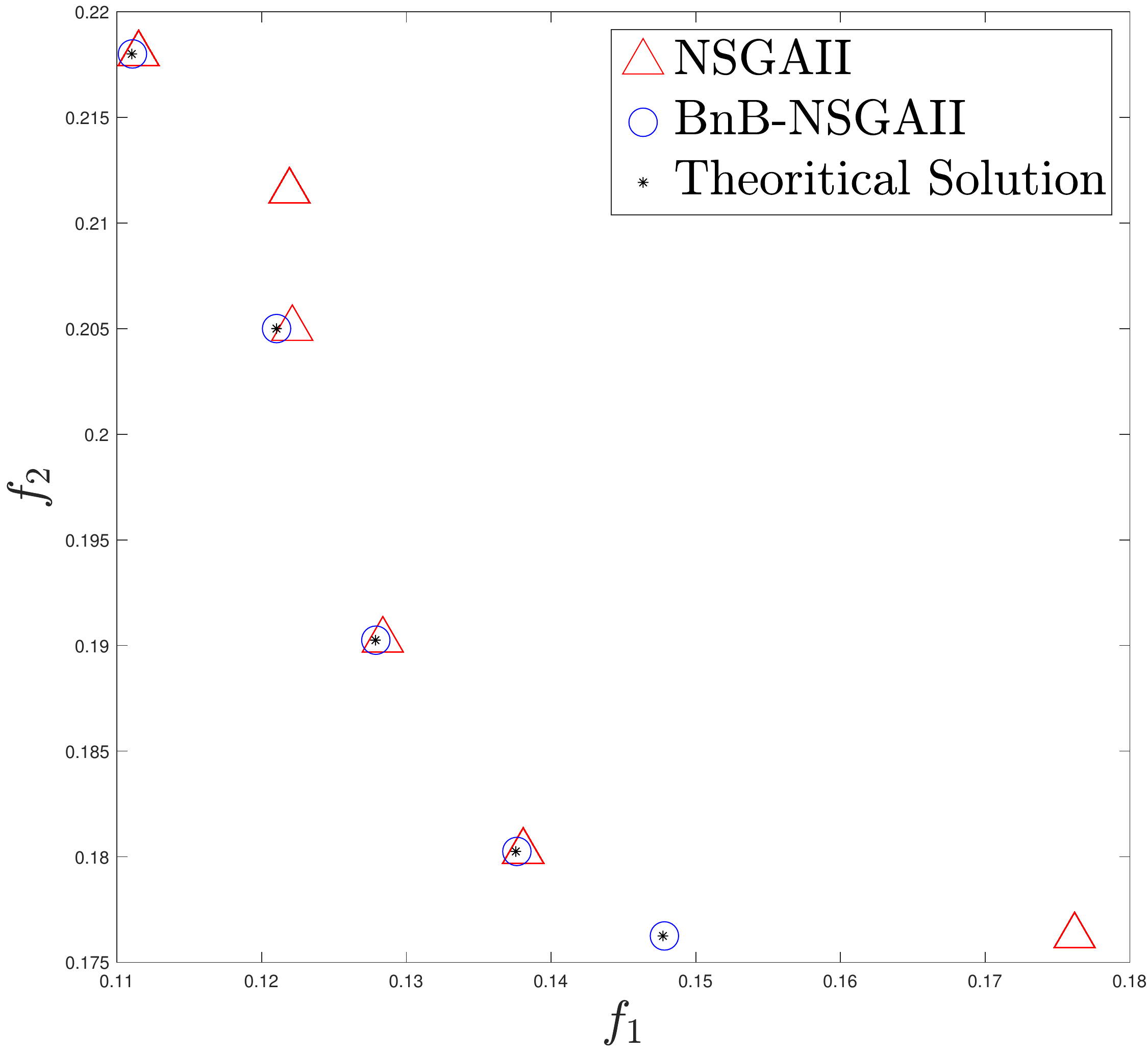}}
	\label{fig:pareto_bearing}
	}
	\subfigure[Coupling problem.]{\resizebox{5cm}{!}{\includegraphics{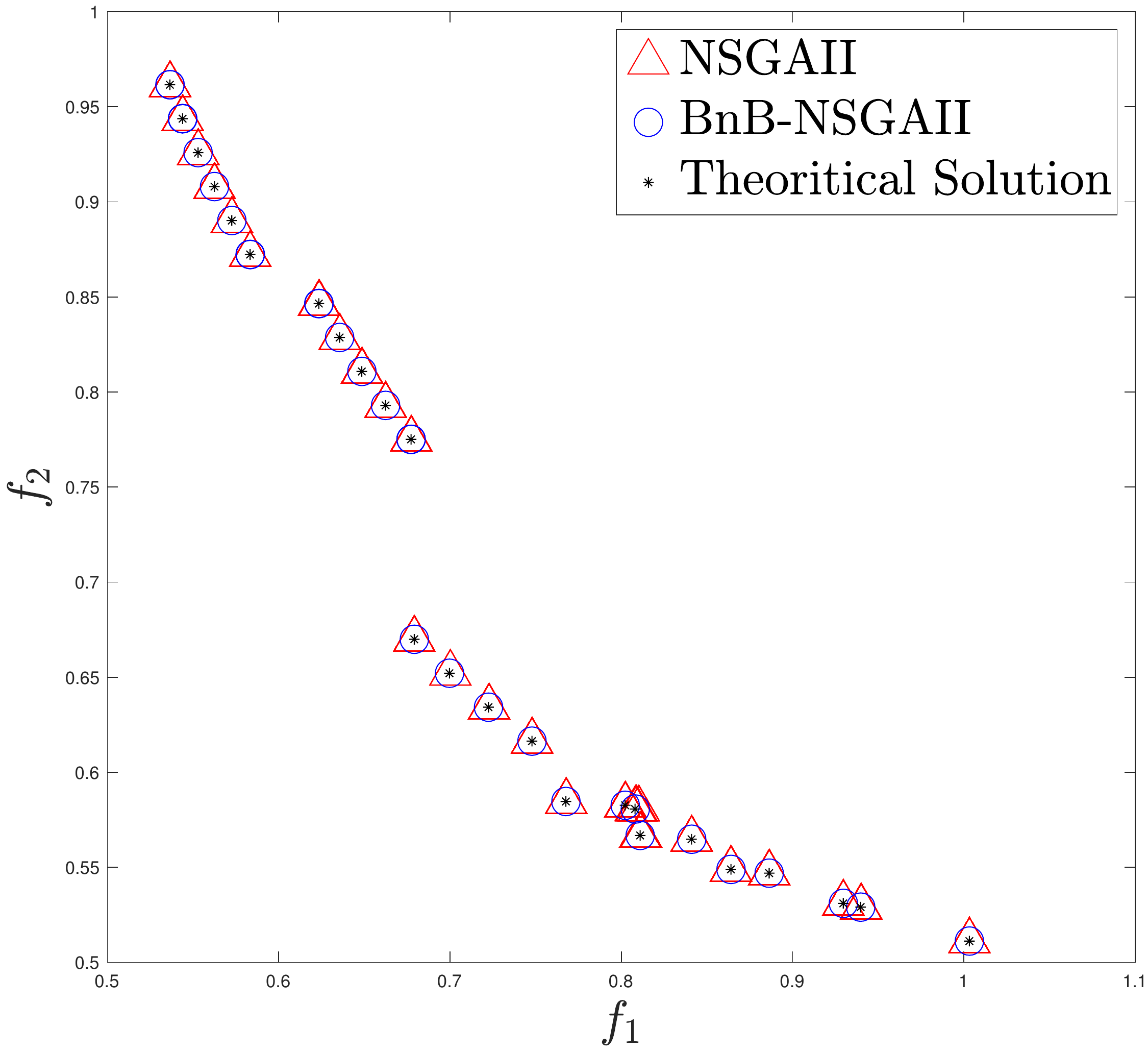}}
	\label{fig:pareto_coupling}
	}\hspace{5pt}
	\subfigure[Brake problem.]{\resizebox{5cm}{!}{\includegraphics{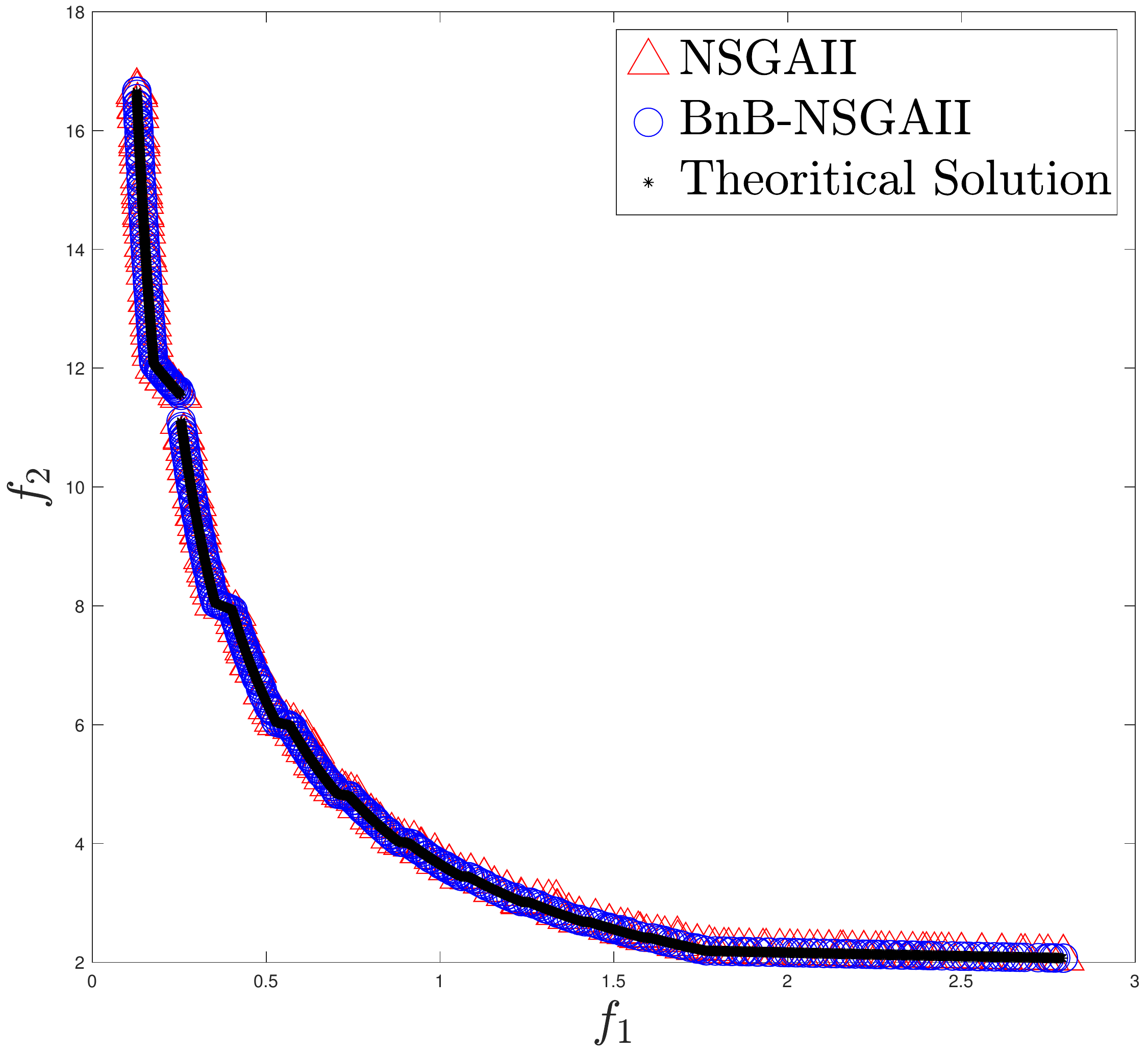}}
	\label{fig:pareto_brake}
	}
	\subfigure[Truss problem.]{\resizebox{5cm}{!}{\includegraphics{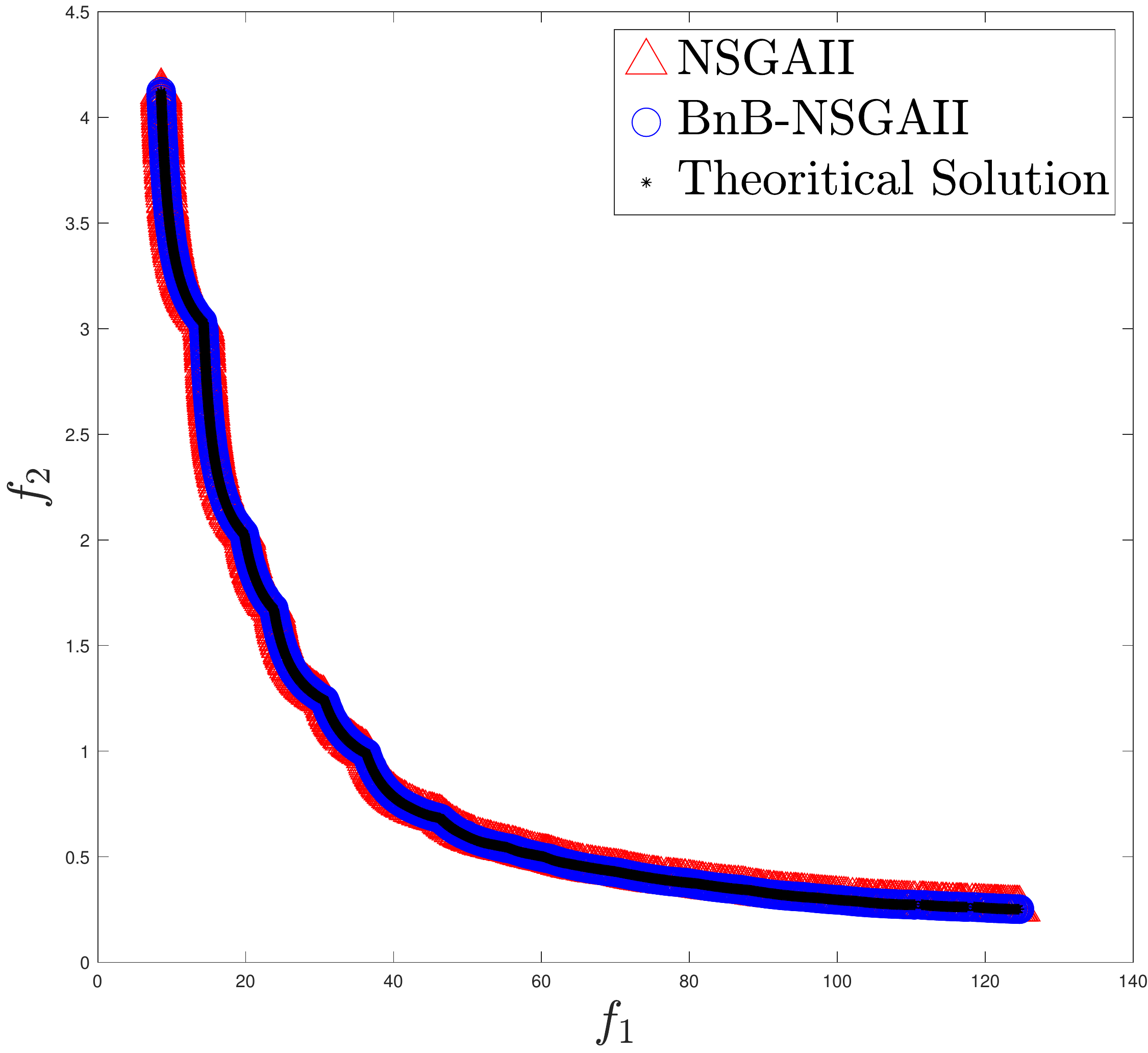}}
	\label{fig:pareto_truss}
	}\hspace{5pt}
	\subfigure[Mela problem.]{\resizebox{5cm}{!}{\includegraphics{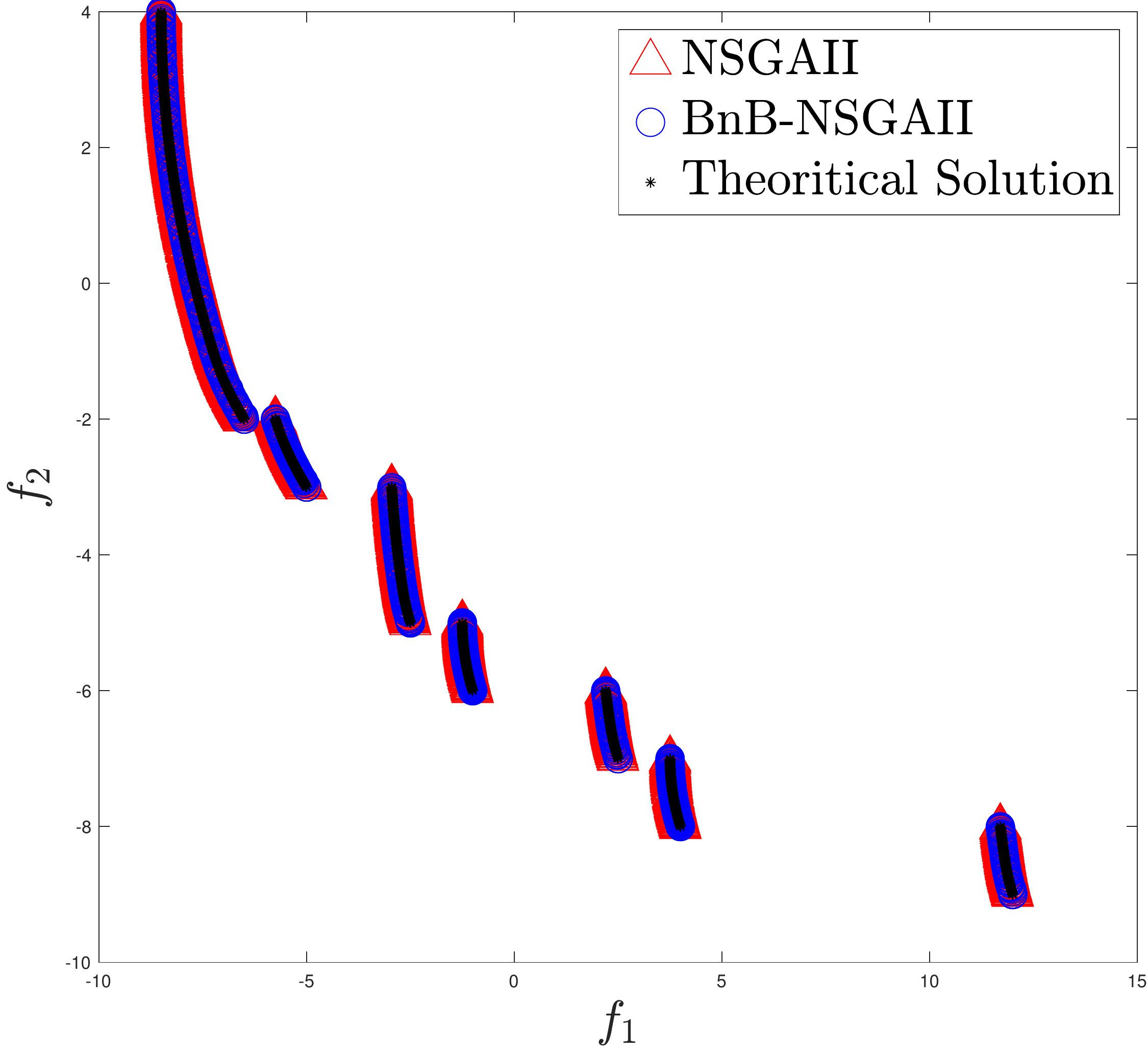}}
	\label{fig:pareto_mela}
	}
	\subfigure[Tong problem.]{\resizebox{5cm}{!}{\includegraphics{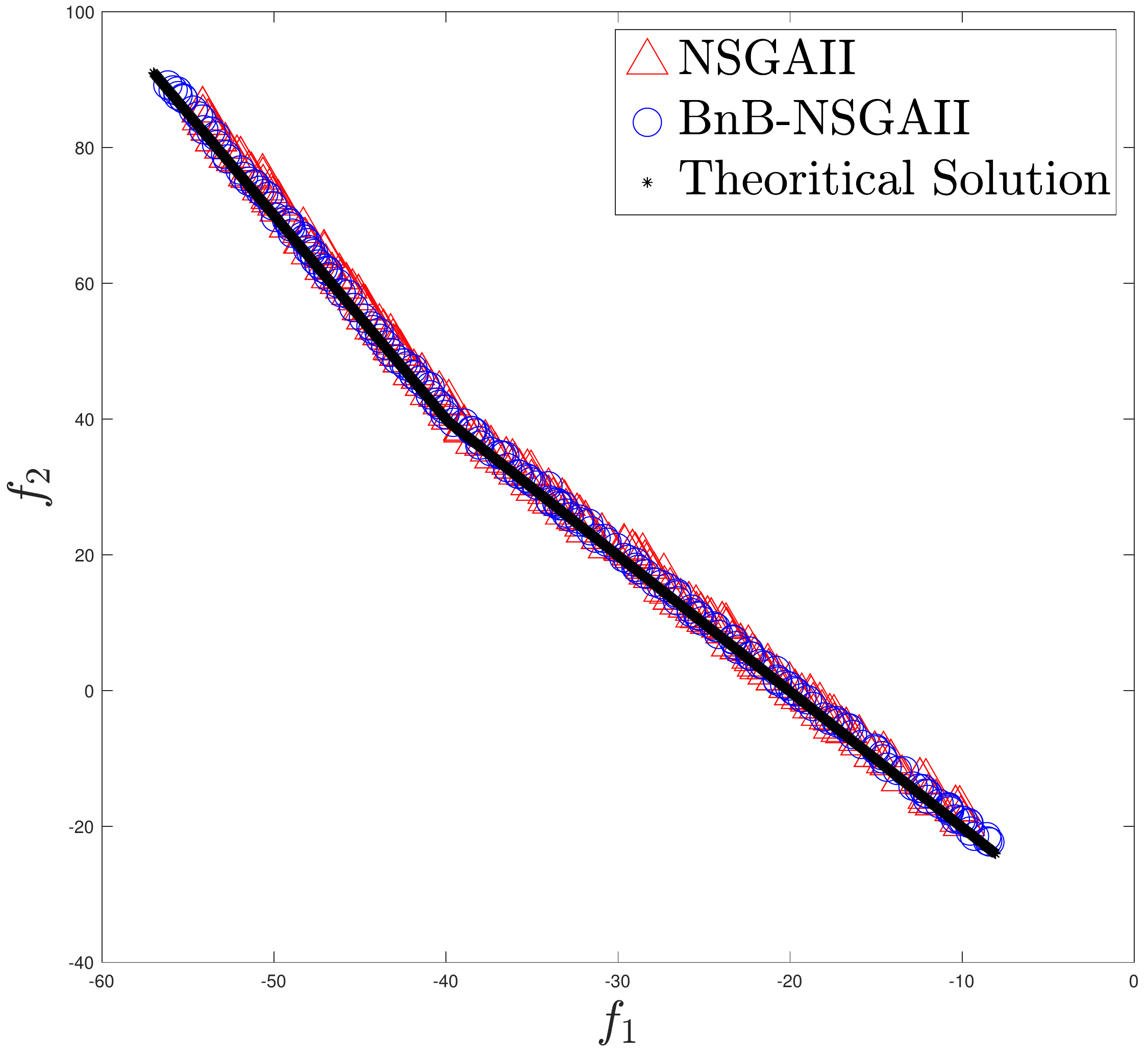}}
	\label{fig:pareto_tong}
	}
	\caption{Pareto front obtained by BnB-NSGAII and NSGAII vs the true Pareto for each problem.}
	\label{fig:pareto}
\end{figure}

 
For continuous/discontinuous Pareto problems, the maximum number of Pareto elements could be obtained by NSGAII is limited by the population size. In Table \ref{tab_results}, it is shown that more elements are found by BnB-NSGAII ($\mid S\mid_\text{BnB-NSGAII}>\mid S\mid_\text{NSGAII}$) since in BnB-NSGAII the number of elements is limited by ($\text{population size} * \text{number of leaf nodes}$).
\begin{table}[h!]
	\tbl{  Mean values of GD, relative spread (d$\Delta$), number of solutions ($\mid S\mid$) and number of evaluations in addition to IR for the selected combinations ID. (Better results are \textbf{emphasized})}{
		\begin{tabular}{llccccccc}
			\toprule
			\textbf{Problem}& \textbf{Method} & \textbf{ID} & \multicolumn{1}{l}{\textbf{GD}} & \multicolumn{1}{l}{\textbf{d$\Delta$}} & \multicolumn{1}{l}{\textbf{$\mid S\mid$}} & \multicolumn{1}{p{1cm}}{\centering  \textbf{Purity}\\ \%} & \multicolumn{1}{p{1cm}}{\centering \textbf{\#Evals}\\ \SI{1e5}{}} & \multicolumn{1}{r}{\textbf{IR}} \\ \midrule
			\multirow{3}{*}{\textit{Gear}} & NSGAII & 1 & \SI{5.6e-2}{} & \SI{1.9e-2}{} & 27   & 97 & 10  & 1        \\ \cmidrule(l){2-9} 
			& \multirow{2}{*}{BnB-NSGAII} & 20 & \SI{3.6e-2}{} & 0       & 28   & 99   & 16  & 24.5 \\ & &
            \textbf{18} & \textbf{0}        & \textbf{0}        & \textbf{28} & \textbf{100} & \textbf{31.2} & \textbf{$\infty$}                             \\ \midrule
			\multirow{4}{*}{\textit{Bearing}}  & NSGAII                      & 1 & \SI{8.8E-04}{} & \SI{6.7E-02}{} & 5    & 92 & 10  & 1        \\ \cmidrule(l){2-9} 
			& \multirow{2}{*}{BnB-NSGAII} & \textbf{19} & \textbf{\SI{1.6E-04}{}} & \textbf{\SI{5.0E-3}{}} & \textbf{5}  & \textbf{100}        & \textbf{5.4} & \textbf{1.59}  \\ & &
22 & \SI{6.2E-5}{} & \SI{2.0E-3}{} & 5    & 100        & 26.6 & 0.82  \\ \midrule
			\multirow{3}{*}{\textit{Coupling}} & NSGAII                      &1          & \SI{2.9E-5} &	\SI{6.0E-1}          & 47           & 34          & 10          & 1                 \\ \cmidrule(l){2-9} 
			& \multirow{2}{*}{BnB-NSGAII} & \textbf{25} & \textbf{\SI{4.2E-5}{}} & \textbf{\SI{4.0E-1}{}} & \textbf{40}  & \textbf{63} & \textbf{9.2}  & \textbf{1.11} \\ & & 
24          & \SI{1.4E-5}
          & \SI{2.8E-1}
          & 35           & 71          & 30.2          & 0.71          \\ \midrule
			\multirow{3}{*}{\textit{Brake}}    & NSGAII                      &2          & \SI{5.6E-4}	&\SI{1.2E-1}
          & 188          & 20          & 10        & 1                 \\ \cmidrule(l){2-9} 
			& \multirow{2}{*}{BnB-NSGAII} & \textbf{12} & \textbf{\SI{2.1E-4}{}} & \textbf{\SI{1.5E-1}{}} & \textbf{399} & \textbf{21}  & \textbf{3.9}  & \textbf{3.68} \\ & & 
6          & \SI{1.4E-4} &	\SI{1.8E-1}
          & 474          & 24          & 22.5          & 0.73         \\ \midrule
			\multirow{3}{*}{\textit{Mela}}     & NSGAII                      & 1          & \SI{3.3E-5}{}	&\SI{8.7E-2}{}
          & 998          & 41          & 10          & 1                 \\ \cmidrule(l){2-9} 
			& \multirow{2}{*}{BnB-NSGAII} & \textbf{10} & \textbf{\SI{8.8E-5}{}
} & \textbf{\SI{8.6E-3}{}
}   & \textbf{835} & \textbf{33} & \textbf{3.6}  & \textbf{5.29} \\ & &
			7          & \SI{2.3E-5}{}	&\SI{6.7E-2}{}
          & 1535         & 40          & 12.1          & 1.13          \\\midrule
			\multirow{2}{*}{\textit{Truss}}    & NSGAII                      & \textbf{1} & \textbf{\SI{9.6E-5}{}
} & \textbf{\SI{1.7E-1}{}
} & \textbf{915} & \textbf{29} & \textbf{10} & \textbf{1}        \\ \cmidrule(l){2-9} 
			& \multirow{1}{*}{BnB-NSGAII} & 7          & \SI{6.4E-5}{}	&\SI{2.4E-1}{}
          & 1438         & 30          & 53.7          & 0.19          \\ \midrule
			\multirow{3}{*}{\textit{Tong}}     & NSGAII                      & 6          & \SI{3.7E-2}{} &	\SI{1.5E-1}{}
          & 141          & 0                 & 10          & 1                 \\ \cmidrule(l){2-9} 
			& \multirow{2}{*}{BnB-NSGAII} & \textbf{10} & \textbf{\SI{7.0E-2}{}
} & \textbf{\SI{2.7E-2}{}
}  & \textbf{119} & \textbf{0}        & \textbf{3.6}  & \textbf{4.85}  \\ & &
			16          & \SI{3.2E-2}{} &	\SI{2.2E-2}{}
          & 181          & 0                 & 36          & 0.78             \\ \bottomrule
		\end{tabular}}
	\label{tab_results}
\end{table}
Figure \ref{fig:pareto} shows the best solution obtained by both methods compared to the true Pareto front for each problem. As shown, BnB-NSGAII is able to approximately converge to the true Pareto front for all the problems. While NSGAII failed to catch all the Pareto elements for gear and bearing problems as shown in Figures \ref{fig:pareto_gear} and \ref{fig:pareto_bearing} respectively. On a detailed examination of the Pareto front of Tong problem Figure \ref{fig:pareto_tong}, it is found that NSGAII does not capture the section of the Pareto front close to the anchor points.

\section{Conclusion}\label{conc}
It was noticed in literature the lack of hybridization of multi-criteria branch and bound (MCBB) with metaheuristics, in contrary to the mono-objective branch and bound. We proposed a general-purpose multi-criteria branch-and-bound based on genetic algorithm (NSGAII) for solving non-convex multi-objective Mixed Integer Non-Linear Programming problems (MINLPs), in particular mechanical design problems. The proposed algorithm and NSGAII were tested on 7 MO-MINLP problems from literature. Results show that the proposed algorithm obtain better solutions for all the problems, although in some cases the computational effort was higher than those for NSGAII. In this aspect, we designed a new metric that relates the quality of the solution to the cost, which we called investment ratio. It was noticed that the results were enhanced by a higher ratio than the computational cost ratio for 6 out of 7 problems, which prove the effectiveness of the proposed method.\par 
\section*{Acknowledgement}
The majority of the numerical experiments were performed on resources provided by the University of Technology of Troyes (UTT).
\section*{Disclosure statement}
No potential conflict of interest was reported by the authors.
\section*{Funding}
This work was supported by the Lebanese University (LU); the regional council of “Grand Est” region, France; and the European Development Fund regional (ERDF).

\bibliographystyle{tfcad}
\bibliography{My_Library}
\appendix
\section{Problems Formulation}\label{prob_form}
\begin{figure}[h!]
\centering
		\includegraphics[width=0.6\linewidth]{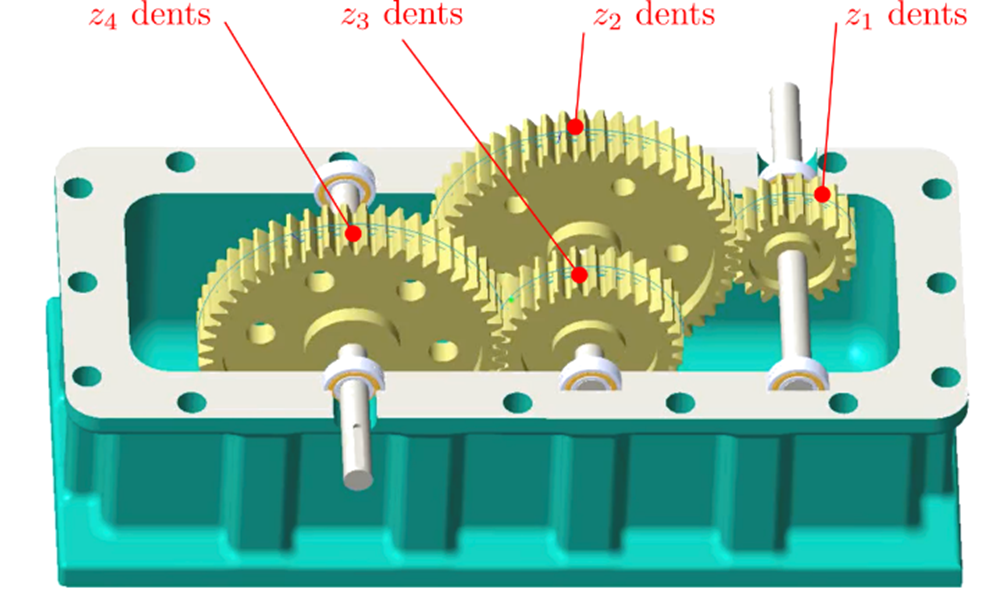}
	\caption{Illustration of gear problem}
	\label{gear_prob}
\end{figure}

The gear train problem is unconstrained. Since the number of teeth must be integers, all 4 variables are strictly integers. The problem is illustrated in figure \ref{gear_prob} and formulated as following:
\begin{mini}<b>
	{}{
		\begin{split}
			f_1 &= 1/6.931 - (z_1*z_2/z_3*z_4)^2\\
			f_2 &= \max\left\lbrace z_1,z_2,z_3,z_4\right\rbrace	
			\end{split}}{}{}
	\addConstraint{}
	\labelOP{pg}
	\addConstraint{12 \leq z_i \leq 60 ;}{}{\quad i=1,2,3,4}
	\addConstraint{ z_i \in  \mathbb{N} ;}{}{\quad i=1,2,3,4}
\end{mini}

\begin{figure}[h!]
	\centering
	\includegraphics[width=0.7\linewidth]{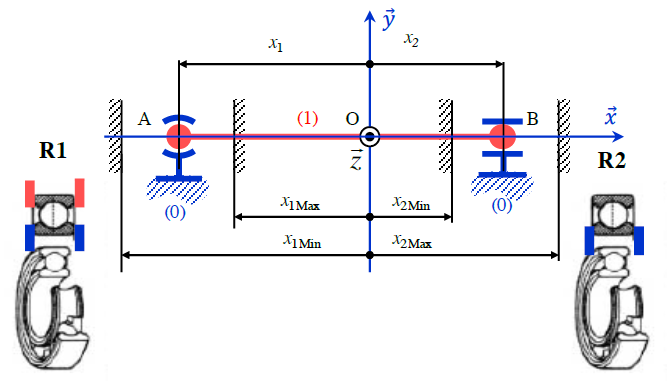}
	\caption{Illustration of bearing problem }
	\label{bearing_prob}
\end{figure}

In the ball bearing pivot link, the ball bearings were chosen from a standardized table of prefabricated sizes. The formulated optimization problem contains 4 variables, 2 continuous/ 2 integer, 12 discrete parameters and 10 inequality constraints. Figure \ref{bearing_prob} illustrates the problem that is formulated as following:\par
 \begin{mini}<b>
 {i_{r1},i_{r2},x_1,x_2}{
		\begin{split}
		F_1&=\frac{\pi\rho_a}{4}\left[\phi_3(i_{r1})^2\phi_5(i_{r1})+\phi_3(i_{r2})^2\phi_5(i_{r2})+\ldots\right.\\
        &\ldots\left(x_1+x_2-\frac{\phi_5(i_{r1})+\phi_5(i_{r2})}{2}\right)\ldots\\
        &\ldots*\max\{\phi_6(i_{r2}),\phi_6(i_{r2})\}^2+\phi_7(i_{r1})+\phi_7(i_{r2}) \\
        F_2&=\phi_1(i_{r1})+\phi_1(i_{r2})
	\end{split}}{}{}
	\addConstraint{}
	\labelOP{pb}
        \addConstraint{C_1=}{\left(XF_{r1}(x_1,x_2)+YF_{a1}\right)\left(\frac{60 L_v\omega}{10^{6}}\right)^\frac{1}{3}-\phi_1(i_{r1})}{\leq 0}
        \addConstraint{C_2=}{F_{r2}(x_1,x_2)\left(\frac{60 L_v\omega}{10^{6}}\right)^\frac{1}{3}-\phi_1(i_{r2})}{\leq 0}
        \addConstraint{C_3=}{ \left(x_{1\text{Min}}+\frac{\phi_5(i_{r1})}{2}\right) -  x_1 }{\leq 0}
        \addConstraint{C_4=}{x_1 - \left(x_{1\text{Max}}-\frac{\phi_5(i_{r1})}{2}\right) }{\leq 0}
        \addConstraint{C_5=}{ \left(x_{2\text{Min}}+\frac{\phi_5(i_{r2})}{2}\right) -  x_2 }{\leq 0}
        \addConstraint{C_6=}{x_2 - \left(x_{2\text{Max}}-\frac{\phi_5(i_{r2})}{2}\right) }{\leq 0}
        \addConstraint{C_7=}{d_{1\text{Min}} - \phi_3(i_{r1}) }{\leq 0}
        \addConstraint{C_8=}{d_{2\text{Min}} - \phi_3(i_{r2}) }{\leq 0}
        \addConstraint{C_9=}{\max\left\lbrace \phi_4(i_{r1}),\phi_4(i_{r2}) \right\rbrace - D_{\text{Max}}}{\leq 0}
        \addConstraint{C_{10}=}{\phi_4(i_{r1}) - \phi_4(i_{r2})}{= 0}
        \end{mini}
        \allowdisplaybreaks
        \begin{align*}
        \text{Where:}&\\
                     &i_{r1},i_{r2} \in \lbrace 1,\ldots,61\rbrace\times\lbrace 1,\ldots,61\rbrace\\\\
        \text{Given :}&\\
                        &\left\lbrace x_{1,2\text{Min}}, x_{1,2\text{Max}}, d_{1,2\text{Min}},D_{\text{Max}},\right.\ldots\\
                        &\left.\ldots \rho_a,\omega, L_v, X_O, Y_O, Z_O, M_O, N_O
                        \right\rbrace
    \end{align*}
With :
    \begin{description}
        \item[$L_{10}$ :] Lifespan expressed in millions of turns [].
        \item[$C$ :] Dynamic load bearing capacity, intrinsic to rolling bearing [\SI{}{\N}].
        \item[$P_{eq}$ :] Radial equivalent load [\SI{}{\N}].
        \item[$F_{a}$ :] Axial load in [\SI{}{\N}].
        \item[$F_{r}$ :] Radial load in [\SI{}{\N}].
        \item[$X$,$Y$ :] A dimensionless coefficient that depends on the bearing type and on the ratio $\frac{F_{a}}{C_{0}}$.
        \item[$C_0$ :] Static load bearing capacity, intrinsic to the rolling bearing in [\SI{}{\N}].
    \item[$B_{1,2}$ :] Bearing width (R1) and (R2) [\SI{}{\mm}].
    \item[$D_{1,2}$ :] Outside diameter of bearings (R1) and (R2) [\SI{}{\mm}].
    \item[$d_{1,2}$ :] Bearing inner diameter (R1) and (R2) [\SI{}{\mm}].
    \item[$x_{1,2\text{Min}}$ :] Minimum limit of bearings positions (R1) and (R2) [\SI{}{\mm}].
    \item[$x_{1,2\text{Max}}$ :] Maximum limit of bearing positions (R1) and (R2) [\SI{}{\mm}].
    \item[$D_{\text{Max}}$ :] Maximum bearing mounting diameter [\SI{}{\mm}].
        \item[$m_{1,2}$ :] Mass of bearings (R1) and (R2) [\SI{}{\gram}].
        \item[$d_{a1,2}$ :] Diameter of the necessary \underline{\textbf{shoulders}} on the shaft for bearings (R1) and (R2) [\SI{}{\mm}].
        \item[$\rho_a$ :] Shaft density [\SI{}{\gram\per\mm^3}]. 
    \end{description}

\begin{figure}[h!]
	\centering
	\includegraphics[width=0.9\linewidth]{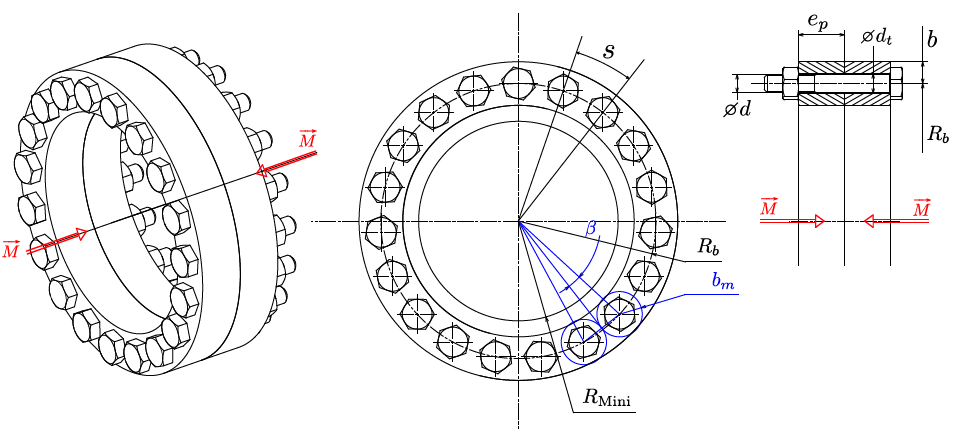}
	\caption{Illustration of coupling problem }
	\label{coupling_prob}
\end{figure}

The coupling with bolted rim contains 1 discrete variable, 1 integer variable, 2 continuous variables, 11 inequality constraints and 5 discrete bolt parameters. Figure \ref{coupling_prob} illustrates the problem that is formulated as following:\par
 \begin{mini}<b>
    {}{
     \begin{split}
     f_1\left(\bm{x}\right)&=\frac{\pi}{2}e_p\left[\rho_J \left(16 R_b\phi_3(i_d)  - N\phi_4(i_d)^2\right) + \rho_b N\phi_1(i_d)^2 \right]\\
        f_2\left(\bm{x}\right)&=\frac{\phi_1(i_d)}{d_{\text{Maxi}}}+\frac{N}{N_{\text{Maxi}}} + \frac{M}{M_{\text{Maxi}}}
     \end{split}
     }{}{}
     \addConstraint{}
     \labelOP{pc}
        \addConstraint{c_1\left(\bm{x}\right)=} {\frac{\alpha_S M}{N R_b K(i_d)} - 1 }{\leq 0}
        \addConstraint{c_2\left(\bm{x}\right)=}{\frac{\phi_3(i_d)N}{k_s \pi R_b} -1 }{\leq 0}
        \addConstraint{c_3\left(\bm{x}\right)=}{\left(R_{\text{Mini}}+\phi_3(i_d)\right) -R_b }{\leq 0}
        \addConstraint{c_4\left(\bm{x}\right)=}{R_b-\left(R_{\text{Maxi}}-\phi_3(i_d)\right)}{\leq 0}
         \end{mini}
         \allowdisplaybreaks
         \begin{align*}
        \text{where:}&\\
                     &\bm{x}=\left\lbrace i_d,N,R_b,M \right\rbrace^T\\
                     &i_d \in \lbrace 1\ldots15 \rbrace\\
                     &N  \in \left[N_{\text{Mini}}, N_{\text{Maxi}}\right]\times\mathbb{N}\\
                     &R_b\in \left[R_{\text{Mini}}, R_{\text{Maxi}} \right]\\
                     &M  \in \left[M_T,M_{\text{Maxi}}\right]\\
        \text{given:}&\\
                        &K(i_d) = \frac{0.9 R_e f_m \pi (\phi_1(i_d) - 0.93815\phi_2(i_d))^2}{4\sqrt{1+3\left(\frac{4\left(0.16\phi_2(i_d)+0.583(\phi_1(i_d) - 0.6495\phi_2(i_d))f_1\right)}{\phi_1(i_d) - 0.93815\phi_2(i_d)}\right)}}\\
                        &k_s =\frac{\sin\left(\frac{\pi}{ N_{\text{Mini}}}\right) - \sin\left(\frac{\pi}{ N_{\text{Maxi}} }\right) }{\frac{\pi}{N_{\text{Mini}}} -\frac{\pi}{N_{\text{Maxi}}}}\\
                        &\left\lbrace N_{\text{Mini}}, N_{\text{Maxi}}, k_s, R_{\text{Mini}} ,R_{\text{Maxi}}, e_p,M_T,\right.\ldots\\
                        &\left.\ldots M_{\text{Maxi}}, \alpha_S, f_1, f_m, R_e, \rho_J, \rho_b \right\rbrace
   \end{align*}
    With:
        \begin{description}
        \item[$C_1$ :] Torsion torque in the bolt body due to tightening [\SI{}{\N\mm}].
        \item[$F_0$ :] Tension force in the bolt body [\SI{}{\N}].
        \item[$d$ :] Nominal thread diameter [\SI{}{\mm}].
        \item[$p$ :] bolt thread pitch [\SI{}{\mm}].       
        \item[$d_2$ :] Average diameter at the thread end of the thread [\SI{}{\mm}].
        \item[$f_1$ :] Coefficient of friction thread [ul].
        \item[$F_{0\text{Maxi}}$ :] Maximum tension in the bolt body [\SI{}{\N}].
        \item[$F_{0\text{Mini}}$ :] Minimum tension in the bolt body [\SI{}{\N}].
        \item[$\alpha_S$ :] Uncertainty coefficient of the clamping tool $>1$ [ul].
        \item[$d_s$ :] Diameter of the "resistant" section of the thread [\SI{}{\mm}].
        \item[$d_3$ :] Diameter of the thread \underline{\textbf{kernel}} [\SI{}{\mm}].
        \item[$C_1$ :] Torsion torque calculated for $F_0=F_{0\text{Maxi}}$.  
        \item[$R_e$ :] Elastic limit of the bolt depending on its quality class [\SI{}{\MPa}].
        \item[$M$ :] Torque transmitted by coupling [\SI{}{\N\m}].
        \item[$N$ :] Number of bolts [ul].
        \item[$R_b$ :] Position radius of $N$ bolts [\SI{}{\mm}].
        \item[$f_m$ :] Coefficient of friction between the plates of the coupling [ul].
        \item[$s$ :] Deviation on the circumference between two bolts [\SI{}{\mm}].
        \item[$b_m$ :] \underline{\textbf{Radial size}} of the clamping tool, depending on the nominal diameter of the bolt $d$ [\SI{}{\mm}].
        \item[$s_m$ :] \underline{\textbf{Circumference}} of the clamping tool, depending on $b_m$ and the installation radius of the $N$ bolt[\SI{}{\mm}].       
        \item[$s_m$ :] Diameter of the holes in the rims, depending on the nominal diameter of the bolt $d$ [\SI{}{\mm}].
        \item[$e_p$ :] Thickness of the coupling rims [\SI{}{\mm}].
        \item[$\rho_J$ :] Density of the rim material [\SI{}{\kg\per\mm^3}].
        \item[$\rho_b$ :] Density of the bolt material [\SI{}{\kg\per\mm^3}].
        \item[$d_{\text{Maxi}}$ :] Diameter of the largest bolts [\SI{}{\mm}].
        \item[$N_{\text{Maxi}}$ :] Maximum number of bolts [ul].
        \item[$M_{\text{Maxi}}$ :] Maximum torque that can be transmitted by adhesion [\SI{}{\N\m}].
    \end{description}
In the disk brake design problem, there are five inequality constraints  that relate to the surface area, length of the brake, pressure,  torque, and temperature. This problem is formulated as following:
\begin{mini}<b>
	{}{
		\begin{split}
			f_1 &= 4.9*10^{-5}(x_1^2-y_1^2)(x_4-1)\\
			f_2 &= \frac{9.82*10^6(x_1^2-y_1^2)}{x_2x_3(x_1^3-y_1^3)}
	\end{split}}{}{}
\labelOP{pbrake}
	\addConstraint{}
	\addConstraint{g_1=}{20-(x_1-y_1)}{\leq 0}
	\addConstraint{g_2=}{2.5(x_3+1)-30}{\leq 0}
	\addConstraint{g_3=}{\frac{x_2}{\pi(x_1^2-y_1^2)}-0.4}{\leq 0}	
	\addConstraint{g_4=}{\frac{2.22*10^{-3}x_2(x_1^3-y_1^3)}{(x_1^2-y_1^2)^2}-1}{\leq 0}	
	\addConstraint{g_5=}{900-\frac{2.66*10^{-2}x_2x_3(x_1^3-y_1^3)}{(x_1^2-y_1^2)}}{\leq0}	
	\addConstraint{}{55 \leq y_1 \leq 80; \quad y_1 \in  \mathbb{N}}{}
	\addConstraint{}{75 \leq x_1 \leq 110 }{}
	\addConstraint{}{1000 \leq x_2 \leq 3000 }{}
	\addConstraint{}{2 \leq x_3 \leq 20 }{}
	\addConstraint{}{x_i \in  \mathbb{R} \quad i=1,2,3}{}
\end{mini} \par
\begin{figure}[h!]
	\centering
	\includegraphics[width=0.7\linewidth]{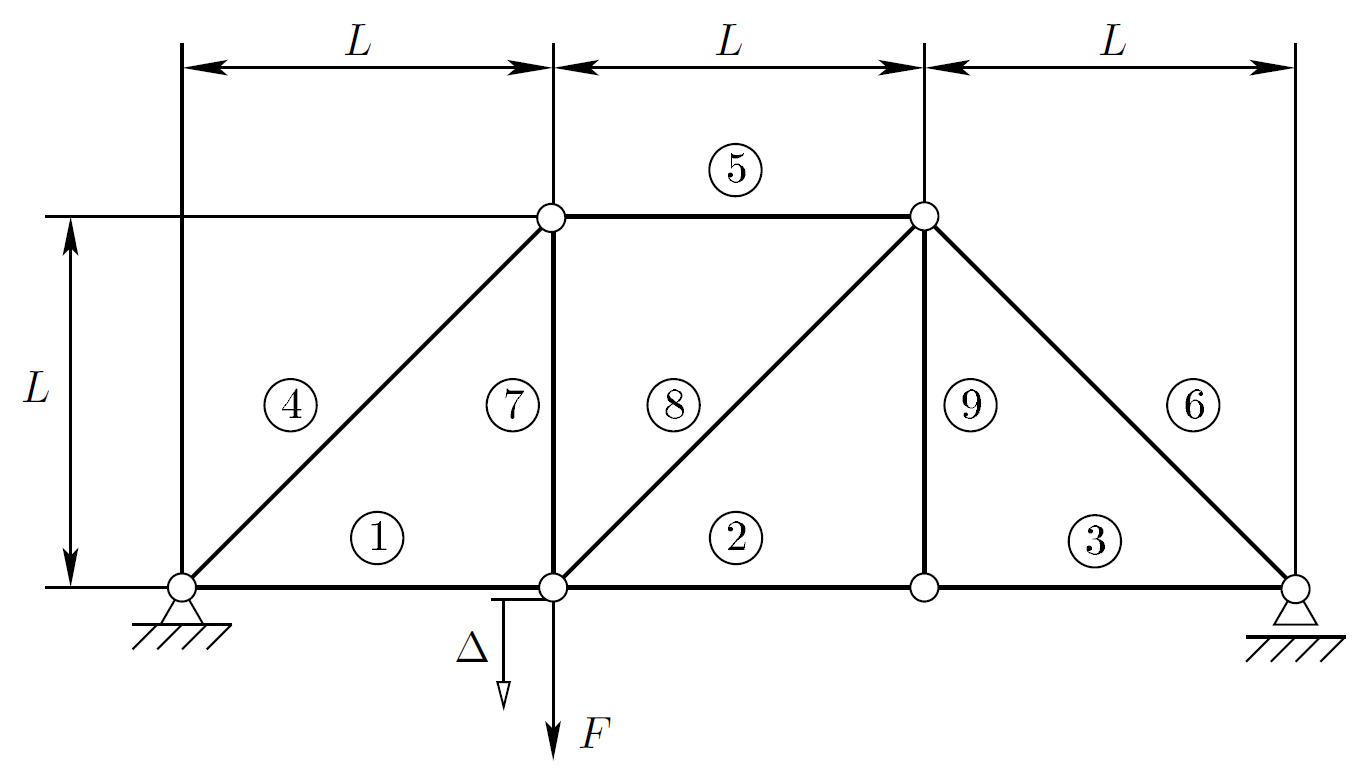}
	\caption{Illustration of truss problem }
	\label{trussprob}
\end{figure}

The nine bar truss problem is unconstrained and contains 3 continuous variables and 9 discrete ones. Figure \ref{trussprob} illustrates the problem that is formulated as following:\par
\begin{mini}<b>
	{}{
		\begin{split}
			f_1 &=(A_1 + A_2 + A_3 + \sqrt{2A_4} + A_5 + \sqrt{2A_6} + A_7 + \sqrt{2A_8} + A_9)\\
			f_2 &= \frac{4}{A_1}+\frac{1}{A_2}+\frac{1}{A_3}+\frac{8\sqrt{2}}{A_4}+\frac{4}{A_5}+\frac{2\sqrt{2}}{A_6}+\frac{4}{A_7}+\frac{2\sqrt{2}}{A_8}
	\end{split}}{}{}
	\addConstraint{}
	\labelOP{ptruss}	
	\addConstraint{}{c_i \leq A_i \leq 10 ;\quad i=1,2,3}
	\addConstraint{}{A_i \in \{1,5,10,15\} ;}{}{\quad i=4,\dots,9}
\end{mini}
where,\[ c_1=\frac{2}{3}, c_2=\frac{1}{3}, c_3=\frac{1}{3}, c_4=\frac{2\sqrt{2}}{3}, c_5=\frac{2}{3}, c_6=\frac{\sqrt{2}}{3}, c_7=\frac{2}{3}, c_8=\frac{\sqrt{2}}{3}, c_9=0\] \par 
Mela problem is unconstrained with 8 binary variables and 2 continuous variables that are bounded. This problem is formulated as following:
\begin{mini}<b>
	{}{
		\begin{split}
			f_1 &=\frac{1}{2}\begin{bmatrix} x & y  \end{bmatrix}^\top\mathbf{G}\begin{bmatrix} x \\ y  \end{bmatrix} +\bm{c}^{1\top}\begin{bmatrix} x \\ y  \end{bmatrix}\\
			f_2 &= \mathbf{c}^{2\top}\begin{bmatrix} x \\ y  \end{bmatrix}
	\end{split}}{}{}
	\addConstraint{}
	\labelOP{pm}	
	\addConstraint{}{-1 \leq x_i \leq 1 ;\quad i=1,2}
	\addConstraint{}{y_j \in \{0,1\} ;}{}{\quad j=1,\dots,8}
\end{mini}
where, $\mathbf{G}, \mathbf{c}^1$ and $\mathbf{c}^2$ are given by:\\
 \[\mathbf{G}= \begin{bmatrix}
1 & -1 &2& 0& 0& 0& 0& 0& 0& 0\\
-1 & 2 & 0& 0& 2& 0& 0& 0& 0& 0\\
0& 0 &3 &0& 2& 0& 0& 0 &0& 0\\
2& 0& 0 &4& 0& 2& 0& 2 &0& 0\\
0& 0& 0 &0 &5 &2 &0& 0& 0& 0\\
0 &0 &0& 0& 0& 6& 0& 0& 0& 0\\
0& 0& 0& 0& 0& 0& 7& 0& 0 &0\\
0 &0 &0 &0& 0& 0& 0& 0& 0& 0\\
0& 0& 0& 0 &0 &2& 0& 0& 0& 0\\
0 &2 &0 &0& 0& 0& 0& 0 &0& 10\end{bmatrix}
\]
\[\bm{c}^1=\begin{bmatrix} -1& -1& 1& -10& 0& 1& -2& 0 &3 &0 \end{bmatrix}^\top\]
\[\bm{c}^2=\begin{bmatrix} 1&  2&  -1&  1&  5&  -2&  0&  6&  0&  -3 \end{bmatrix}^\top\]\par
In Tong problem, it involves 3 continuous variables, 3 binary ones and 9 constraints. The problem is formulated as following:
\begin{mini}<b>
	{}{
		\begin{split}
			f_1 &= x_1^2-x_2+x_3+3y_1+2y_2+y_3\\
			f_2 &= 2x_1^2+x_2-3x_3-2y_1+y_2-2y_3
	\end{split}}{}{}
	\labelOP{pt}
	\addConstraint{}
	\addConstraint{g_1=}{3x_1-x_2+x_3+2y_1)}{\leq 0}
	\addConstraint{g_2=}{4x_1^2+2x_1+x_29x_3+y_1+7y_2}{\leq 40}
	\addConstraint{g_3=}{-x_1-2x_2+3x_3+7y_3}{\leq 0}	
	\addConstraint{g_4=}{-x_1+12y_1}{\leq 10}	
	\addConstraint{g_5=}{x_1-2y_1}{\leq 5}
	\addConstraint{g_6=}{-x_2+y_2}{\leq 20}
	\addConstraint{g_7=}{x_2-y_2}{\leq 40}	
	\addConstraint{g_8=}{-x_3+y_3}{\leq 17}	
	\addConstraint{g_9=}{x3-y_3}{\leq 25}	
	\addConstraint{}{-100 \leq x_i \leq 100;\; i=1,2,3;\;}{x_i \in  \mathbb{R}}
	\addConstraint{}{y_j \in  \{0,1\}; \; j=1,2,3}{}
\end{mini} \par
\section{Determination of true Pareto}\label{true_par}
The true Pareto solution is obtained in the following manner, first, the problem is solved by discretization the continuous variables, then the combinations of integer variables that gives part of the Pareto are recorded. By fixing those combinations iteratively, the problem is solved exactly using the scalarization method as \ProbOpt{MO-NLP}. Cubic interpolation is then used to obtain a uniform discretization step for each continuous part of the Pareto front. Finally, Pareto refining is applied to remove all dominated solutions.\par 
For the convergence metric GD in the case of continuous/discontinuous Pareto, the value of $d_i$ is affected by the discretization step used to generate the true Pareto front. Thus although the tested solver might find a solution that lies on the theoretical front, the distance $d_i$ might be greater than 0 because of the discretization step. To avoid this, first we chose $\epsilon=10^{-4}$ as discretization step, then for each $n$ non-dominated solution $d_n\leq\epsilon$ since $n$ must be between ideal and nadir points. Thus we can assume that $d_n=0$ for a more reasonable approximation of GD. However this assumption might lead to inaccurate GD value, but since we are using this assumption on both compared solvers, we can consider it impartial. For the relative spread metric, to minimize $\Delta_\text{True}$, the true Pareto is discretized uniformly, i.e. the euclidean distance between every two consecutive solutions laying on a continuous section of the front is constant.\par 
\end{document}